\documentclass[a4paper]{article}
\usepackage{inputenc}

\usepackage{amsmath,amsfonts,amssymb}
\usepackage{graphicx}
\usepackage{color}
\usepackage{caption} 
\usepackage{changes}
\usepackage{cancel}
\usepackage{ulem}
\usepackage{lineno}
\usepackage{tabularx}
\usepackage{jabbrv}

\usepackage{mathtools}
\mathtoolsset{showonlyrefs}

\usepackage[textsize=footnotesize,backgroundcolor=yellow!70,bordercolor=orange]{todonotes}
\definecolor{refkey}{rgb}{.35,.75,0}
\definecolor{labelkey}{rgb}{.15,.55,0}

\def\R{{\mathbb R}}

\def\L{{\mathcal L}}

\def\I{{\mathcal I}}

\let\theta\vartheta

\def\vec#1{{\mathbf{#1}}}
\def\vecg#1{{{#1}}}

\def\wh{\vecg{w}_h}

\DeclareMathOperator{\sgn}{sgn}

\title{A level-set multigrid technique for nonlinear diffusion in the numerical simulation of marble degradation under chemical pollutants}
\author{A.~Coco \thanks{(Corresponding author) School of Engineering, Computing and Mathematics, Oxford Brookes University, Wheatley campus, Oxford}, 
M.~Semplice\thanks{Dipartimento di Matematica, Universit\`a di Torino, via C.~Alberto, 10, Torino}
\thanks{Centro Speciale di Scienze e Simbolica dei Beni Culturali, Universit\`a dell'Insubria, Chiostro di S.~Abbondio, via Regina Teodolinda, Como} 
and 
S.~Serra Capizzano\thanks{Department of Humanities and Innovation, Universit\`a dell'Insubria, via Bossi 5, 22100 Como}
\thanks{Centro Speciale di Scienze e Simbolica dei Beni Culturali, Universit\`a dell'Insubria, Chiostro di S.~Abbondio, via Regina Teodolinda, Como}
}

\begin{document}
\maketitle

\begin{abstract}
Having in mind the modelling of marble degradation under chemical pollutants, e.g.~the sulfation process, we consider governing nonlinear diffusion equations and their numerical approximation.
The space domain of a computation is the pristine marble object. In order to accurately discretize it while maintaining the simplicity of finite difference discretizations, the domain is described using a level-set technique. A uniform Cartesian grid is laid over a box containing the domain, but the solution is defined and updated only in the grid nodes that lie inside the domain, the level-set being employed to select them and to impose accurately the boundary conditions.
We use a Crank-Nicolson scheme in time, while for the space variables the discretization is performed by a standard Finite-Difference scheme for grid points inside the domain and by a ghost-cell technique on the ghost points (by using boundary conditions).
The solution of the large nonlinear system is obtained by a Newton-Raphson procedure and a tailored multigrid technique is developed for the inner linear solvers.
The numerical results, which are very satisfactory in terms of reconstruction quality and of computational efficiency, are presented and discussed at the end of the paper.
\end{abstract}

\section{Introduction}
Quantitative forecasts of damage by gaseous pollutants to monuments are becoming more and more important, since they allow to schedule monitoring, preservation and, when needed, restoration activities in the management of cultural heritage \cite{dTPSFCFdM:16:airpollutionimpact}. 
Phenomena involving chemical reactions of the constitutive material of a work of art with chemicals in the surrounding environment have long been recognized very important for the damage to cultural heritage sites. In a recent review of the related mathematical models \cite{reviewSQB:18}, the employment of models based on partial differential equations has been advocated for the next generation of models with regulatory powers.

Differential models in this field typically involve a, possibly nonlinear, diffusion term describing the penetration of the gas in the bulk material, coupled with reaction terms modelling the chemistry of the interaction. A typical example is the sulfation process that turns marble into gypsum, for which a model was proposed in \cite{ADN:sulfation}. More recent models include the effects of the Darcy velocity \cite{AFNT:07:darcymodel}, of the surface rugosity \cite{BCFGN:19}; a kinetic approach to the derivations of the models of sulfation was proposed in \cite{ABST:12:kinetic}.
More complex models including free boundaries can take into account the swelling of the material \cite{CFN:08:swelling} or treat appropriately the heterogeneity of the crust layer \cite{Nikolopoulos:14:mushy}. Modelling of layered material with moving interfaces has been exploited also for copper corrosion in \cite{CdFN:14:copper}.

In this paper we focus on the model of \cite{ADN:sulfation}, since it is quite simple but yet contains the more relevant numerical difficulties.
Previous numerical work on this model have considered the one-dimensional version of the model \cite{GSNF:08:num&lab}
or Cartesian grids in two space dimensions
\cite{Matteo:monum}.
Of course this kind of meshes cannot stand one of the main difficulties of this kind of computations, which is the accurate discretization of the domain. 
This is quite relevant, since, despite the simplicity of the numerical techniques in \cite{Matteo:monum}, the computational model have shown the importance of two-dimensional effects near corners and other sharp features of the domain.
In a real case $\Omega$ should coincide with the pristine work of art and only rarely this can be accurately represented on a Cartesian grid.

One of the simplest methods to overcome this difficulty consists of approximating $\Omega$ by small cuboids, whose size and shape are eventually adapted close to the boundary of $\Omega$ in such as way that the most external corners lie on the boundary. This approach is adopted by the Shortley-Weller discretization~\cite{Shortley-Weller:discretization}, providing the simplest approach falling under the class of boundary-fitted grid methods, where the grid is suitably adapted to the boundary of the domain.
Although the Shortley-Weller discretization is designed for Dirichlet boundary conditions and can be highly accurate for the solution and its gradient~\cite{seo2018convergence}, the extension to Neumann boundary conditions is not straightforward.	

More accurate and well-known boundary-fitted methods are represented by the Finite Element Methods (FEM), successfully adopted in several scientific contexts (e.g., \cite{brenner2007mathematical, brezzi2012mixed, girault2012finite, Babuska:FEM, Bramble:FEM, huang2002mortar, hansbo2004finite, dryja2005neumann, gross2007extended, dolbow2009efficient}). 
However, fitting the mesh to a complex domain with several corners or highly variable curvature might be computationally demanding.
Another alternative for working with a Cartesian structure is furnished by the Isogeometric approach, adopted for example in~\cite{xu2011parameterization, hesch2012isogeometric}, where an evident difficulty is given by the need of using several patches when the domain is complicate and of course this is a concrete possibility when treating the degradation e.g. of a statue from our cultural heritage.

For all these reasons, the mathematical models proposed in this paper, where the domain coincides with the realistic monument with its sharp features, would be more efficiently solved by numerical approaches where the boundary is embedded in a steady Cartesian grid and implicitly described by a level-set function. Another advantage of this approach is that it would be easier to generalize the methods to the case of evolving boundaries or the presence of internal interfaces, advocated by the more modern models in \cite{CFN:08:swelling,Nikolopoulos:14:mushy,CdFN:14:copper}.

The first methods falling under this category were the Immersed Boundary Method~\cite{Peskin:IBM} and the Immersed Interface Methods~\cite{LeVequeLi:IIM}, proposed to model blood flows in the heart.

More recent numerical approaches to discretise partial differential equations on complex domains in a Cartesian grid are the Ghost-Fluid Methods proposed in~\cite{Fedkiw:GFM, Gibou:Ghost, Gibou:fourth_order, Gibou:fluid_solid}, where the solution and the Dirichlet boundary condition are extrapolated to define ghost values outside the domain in order to maintain a standard discretization stencil on internal grid points without compromising the overall accuracy order.

An improved version that accounts for Neumann boundary conditions was proposed by Coco and Russo in~\cite{CocoRusso:Elliptic, coco2018second, coco2012second}. In this method a high order accuracy is achieved not only for the solution but also for the gradient of the solution. 

Other recent methods for Neumann boundary conditions were proposed in~\cite{helgadottir2015imposing}, while the accuracy order for the gradient of the solution has been improved in~\cite{bochkov2019solving}.

Among the other unfitted-boundary approaches, we mention
the matched interface (MIB) method~\cite{Zhou:MIB},
the Immersed Finite Volume Methods (IFVM)~\cite{Oevermann:FV_sharpInterfaceFV, ewing1999immersed},
the arbitrary Lagrangian Eulerian method (ALE)~\cite{FormaggiaNobile:ALE, Donea:ALE},
and the penalization methods~\cite{lacanette2009eulerian, Angot:penalization, Iollo:penalization}.

Within this paper we consider a novel numerical technique for the approximation of nonlinear (possibly degenerate) parabolic equations, which relies on the finite difference discretization and efficient solvers of \cite{Matteo:monum,DSS,DSS-2} and on the level-set domain description and handling of boundary conditions of \cite{CocoRusso:Elliptic}.
As in \cite{Matteo:monum}, the time discretization is the implicit Crank-Nicolson, a large nonlinear system at each time step is solved by a Newton-Raphson procedure, with a tailored multigrid technique for the linear systems.
The spatial discretization is achieved by finite differences on a uniform Cartesian grid and, in the bulk of the domain, the numerical scheme is the same as in \cite{DSS}. However, here, the domain can be of arbitrary shape and is implicitly defined by $\Omega=\{\vec{x} \text{, s.t. } \varphi(\vec{x})<0\}$, where the level-set function $\varphi:\R^N\to\R$ is known at least at the grid nodes.
The grid nodes are defined, according to $\varphi$, as internal (those inside $\Omega$), ghosts (first layer of points around the internal ones) and external.
The method of \cite{DSS} is applied only on the internal grid points. In order to close the method, the resulting nonlinear system of equations is augmented, as in \cite{CocoRusso:Elliptic}, by the equations expressing the fulfillment of the boundary conditions on $\partial\Omega$ in terms of the ghost values and of their first internal neighbour points. The resulting system is then solved by Newton-Raphson and the special smoothing technique of \cite{CocoRusso:Elliptic} is employed in the multigrid linear solver.

The outline of the paper is the following.
In \S\ref{sec:model} we introduce the mathematical model.
The numerical method is discusses in \S\ref{sec:method}, discussing the details of the time discretization in \S\ref{ssec:method:time}, the space discretization in \S\ref{ssec:method:levelset} and \S\ref{ssec:method:space}, the Newton-Raphson solver in \S\ref{ssec:method:newton} and the multigrid method in \S\ref{ssec:method:multigrid}. The numerical tests of \S\ref{sec:numtest} include  accuracy and efficiency tests, as well as examples of application to nontrivial geometries \S\ref{ssec:numtest:realgeometry}. Finally, the main conclusions of the paper and perspectives for future work are discussed in \S\ref{sec:conclusion}.

The numerical results, which are very satisfactory both from the viewpoint of the reconstruction quality and of the computational efficiency, are presented and discussed at the end of the paper.

\section{Mathematical Model} \label{sec:model}
Here we recall briefly the model of marble sulfation introduced in \cite{ADN:sulfation}, referring the reader to the
original paper for the details and more comprehensive study of the
properties of the solutions. In \cite{ADN:sulfation}, the authors
consider the (simplified) chemical reaction
\begin{equation*}\label{eq:REAC}
\mathrm{CaCO_3} + \mathrm{SO_2} +\frac12\mathrm{O_2} +2\mathrm{H_2O}
\longrightarrow \mathrm{CaSO_4}\cdot2\mathrm{H_2O} + \mathrm{CO_2}.
\end{equation*}
to account for the transformation of $\mathrm{CaCO_3}$ of the marble
stone into $\mathrm{CaSO_4}\cdot2\mathrm{H_2O}$ (gypsum), that is
triggered in a moist atmosphere by the availability of $\mathrm{SO_2}$
at the marble surface and inside the pores of the stone.

Letting $\Omega\subset\R^d$ represent the pristine marble piece, the equations governing the process of marble sulfation are:
\begin{equation}\label{maineq}
\left\{
\begin{array}{rcll}
\displaystyle \frac{\partial \left( \phi(c)s \right)}{\partial t} &=& \displaystyle -\frac{a}{m_c} \phi(c) \, s \, c+ d\, \nabla \cdot (\phi(c) \nabla s) & \mbox{ in } \Omega \times [0,T]\\
\displaystyle \frac{\partial c}{\partial t} &=& \displaystyle -\frac{a}{m_s} \phi (c) \, s\, c & \mbox{ in } \Omega \times [0,T] \\
s(\vec{x},t) &=& s_b & \mbox{ for } (\vec{x},t) \in \partial \Omega \times [0,T] \\
s(\vec{x},0) &=& s_0(\vec{x}) & \mbox{ for } \vec{x} \in \Omega \\
c(\vec{x},0) &=& c_0(\vec{x}) & \mbox{ for } \vec{x} \in \Omega .
\end{array}
\right.
\end{equation}

where $c(t,\vec{x})$ represents the marble concentration, initially set to $c_0$, and $s(t,\vec{x})$ is the gas concentration, initially set to $s_0$. The evolution equations \eqref{maineq} are characterized by a diffusion term for $s$, that is nonlinear since the diffusion coefficient depends on $c$, and by reaction terms coupling the two variables. Boundary conditions are of Dirichlet type and impose a value of $s$ on the boundary of $\Omega$, representing the pollution level of the surrounding air.
It is assumed that, the gypsum concentration is $1-c(t,\vec{x})$ and so, as time goes by, the calcium carbonate concentration is reduced from the initial value $c_0$, as $\mathrm{CaCO_3}$ is progressively replaced by a gypsum crust that forms on the outer shell of the monument and whose thickness and evolution are of interest to the managers of cultural heritage.

The time evolution is described by the diffusion term in the gas equation and by the reaction terms in both differential equations. The porosity $\phi$ of the material controls the diffusion of the gas in the pores of the marble. Since marble and gypsum have different porosities, $\phi$ is not a constant but a function of the $c(\vec{x},t)$, making the diffusion equation a nonlinear one. For simplicity, as in \cite{ADN:sulfation},
we assume that
\[
\phi(c(\vec{x},t)) = \alpha \, c(\vec{x},t) + \beta, \quad \mbox{ with } \alpha = 0.01 \mbox{ and } \beta = 0.1.
\]
We point out that more complex relations may be employed, since our method does not rely strongly on the linearity of the above relation.

The third equation is the boundary condition for $s$ and describes the condition surrounding the work of art. 
We assume for simplicity a Dirichlet boundary condition, although the method proposed in this paper can be extended to more sophisticated boundary conditions on fluxes such as in the free boundary model proposed in 1D in~\cite{CFN:08:swelling}.
We observe that no boundary condition is needed for $c$, since the second equation of \eqref{maineq} does not involve spatial derivatives. Last two equations represent the initial conditions.

Although the model can be cast in higher dimensions, in this paper we focus on the 2D case for simplicity, since it contains already the main difficulties of the full 3D model.
Let $D = [-L,L]^2$ be the computational domain, $\Omega \subset D$ the domain representing the marble monument, $\Gamma = \partial \Omega$ the boundary of the domain.
The numerical method will consider a regular grid in $D$ and a level-set function defined on the grid will be used both to detect the grid points inside $\Omega$  and to describe the exact location and outward normal for the boundary of $\Omega$.

\section{Numerical method}\label{sec:method}
\subsection{Time discretization}\label{ssec:method:time}
Equations \eqref{maineq} are discretized in time using the second order accurate Crank-Nicolson scheme
\begin{equation}\label{discTime}
\left\{
\begin{array}{rcll}
\displaystyle \frac{ \phi(c^{(n+1)})s^{(n+1)} - \phi(c^{(n)})s^{(n)}}{\Delta t} &=& \displaystyle\frac{\mathcal{L}^s \left(s^{(n)},c^{(n)}\right) + \mathcal{L}^s \left(s^{(n+1)},c^{(n+1)}\right)}{2} & \mbox{ in } \Omega \\
\displaystyle \frac{c^{(n+1)}-c^{(n)}}{\Delta t} &=& \displaystyle\frac{\mathcal{L}^c \left(s^{(n)},c^{(n)}\right) + \mathcal{L}^c \left(s^{(n+1)},c^{(n+1)}\right)}{2} & \mbox{ in } \Omega \\
s^{(n+1)} &=& s_b & \mbox{ on } \partial \Omega \\
s^{(0)} &=& s_0 & \\
c^{(0)} &=& c_0 & 
\end{array}
\right.
\end{equation}
where
\[
\L^s (s,c) = \displaystyle -\frac{a}{m_c} \phi(c) \, s \, c+ d\, \nabla \cdot (\phi(c) \nabla s), \quad
\L^c (s,c) = \displaystyle -\frac{a}{m_s} \phi (c) \, s\, c
\]
are the differential operators representing the right-hand side of \eqref{maineq},
while $s^{(n)}=s^{(n)}(\vec{x})$ and $c^{(n)}=c^{(n)}(\vec{x})$ are functions of space only and represent approximations of the solutions $s(\vec{x},t)$ and $c(\vec{x},t)$, respectively, at time $t_n= n \, \Delta t$, for $n=0,\ldots, N_t$ (with $\Delta t = T/N_t$).

\subsection{Level-set function}\label{ssec:method:levelset}
The domain $\Omega$ is represented by an auxiliary level-set function $\varphi \colon D \to \R$:
\[
\Omega = \left\{ \vec{x} \in \R^2 \colon \varphi(\vec{x}) < 0 \right\}, \quad \Gamma = \partial \Omega = \left\{ \vec{x} \in \R^2 \colon \varphi(\vec{x}) = 0 \right\}.
\]
Level-set methods have been firstly proposed by Sethian and Osher~\cite{osher1988fronts} and widely used since then to solve PDEs in complex-shaped geometries and moving domains~\cite{osher2006level}.
The level-set function is convenient to gather geometric information, such as the outward unit normal vector $\vec{n}$ and the curvature $\kappa$ of the boundary $\Gamma$:
\begin{equation}\label{ncurv}
\vec{n} = \frac{\nabla \varphi}{|\nabla \varphi|}, \quad \kappa = \nabla \cdot \vec{n}.
\end{equation}
Among the infinite level-set functions that can describe a domain $\Omega$, the signed distance function
$\phi_d(\vec{x}) = \pm \text{dist}(\vec{x},\Gamma)$ has the advantage to maintain a stable algorithm since steep or shallow gradients are avoided. If the signed distance function $\phi_d$ is not available, it can be approximated from a generic level-set function $\phi$
by approximating the solution of the following PDE up to steady state (reinitialization step):
\begin{equation*}
\frac{\partial \phi_d}{\partial \mu} = \sgn(\phi) \left(1- \left| \nabla \phi_d \right| \right), \quad \phi_d=\phi \; \mbox{ when } \; \mu=0.
\end{equation*}
where $\mu$ is a fictitious time parameter~\cite{osher2006level}.
The numerical methods proposed in this paper are designed for generic level-set function $\phi$ and therefore do not rely on the availability of the signed-distance function.

We observe that in real applications the level-set function must be derived from image data, usually starting from a cloud of points given by a 3D laser scanner. The level-set function can be obtained as follows. First, a provisional level-set function is defined by setting $\phi=-1$ inside the domain and $\phi = 1$ outside the domain, introducing a stepwise effects on the boundary. 
Then, the level-set function is slightly diffused (few iterations of diffusion equation), and the stepwise effects on the boundary disappear without affecting the important features of the monument.
In this paper, we adopt the simplified 2D version of this technique to pixelated black and white images representing a woman head profile, a shark and a necklace (see \S\ref{ssec:numtest:realgeometry}).

\subsection{Spatial discretization}\label{ssec:method:space}
Discretization in space is performed by a standard finite-difference scheme for the grid points lying inside the domain $\Omega$ and by a ghost-cell technique on the ghost points (to impose high order accurate boundary conditions). 

Let $D=[-L,L]^2$ be the computational domain and let $N\geq 1$ be the number of intervals in each direction. We call $h=2L/N$ the spatial step. Observe that we are assuming for simplicity that $\Delta x = \Delta y = h$, although the method can be easily generalised to the case $\Delta x \neq \Delta y$.

The set of grid points is $D_h = \left\{ (x_i,y_j) \in \R^2 \colon x_i =-L+ i\, h, y_j= -L+j\, h, \mbox{ for } i,j=0, \ldots, N \right\}$. Let $\Omega_h = \Omega \cap D_h$ be the discrete counterpart of $\Omega$. 

\begin{figure}[!hbt]
 \begin{minipage}[c]{0.45\textwidth}
   	\centering
   	\captionsetup{width=0.80\textwidth}
		\includegraphics[width=0.99\textwidth]{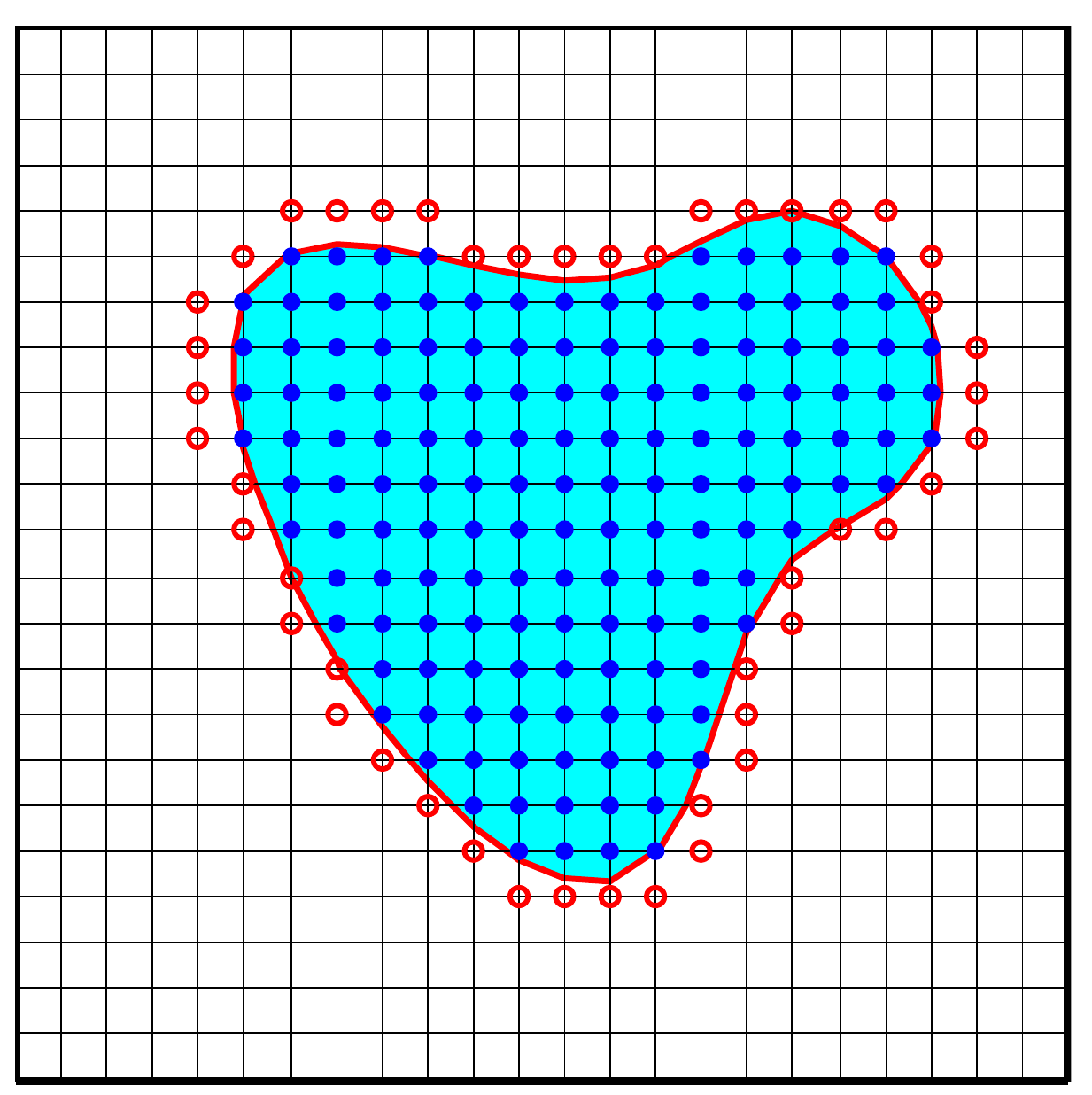}
		\caption{ \footnotesize{Internal grid points $\Omega_h$ (blu filled circle) and ghost points $\Gamma_h$ (red empty circle), according to the definition of ghost points \eqref{ghostdef}. } }
	\label{fig:domain}
 \end{minipage}
 \ \hspace{2mm} \hspace{3mm} \
 \begin{minipage}[c]{0.45\textwidth}
   	\centering
   	\captionsetup{width=0.80\textwidth}
		\includegraphics[width=0.99\textwidth]{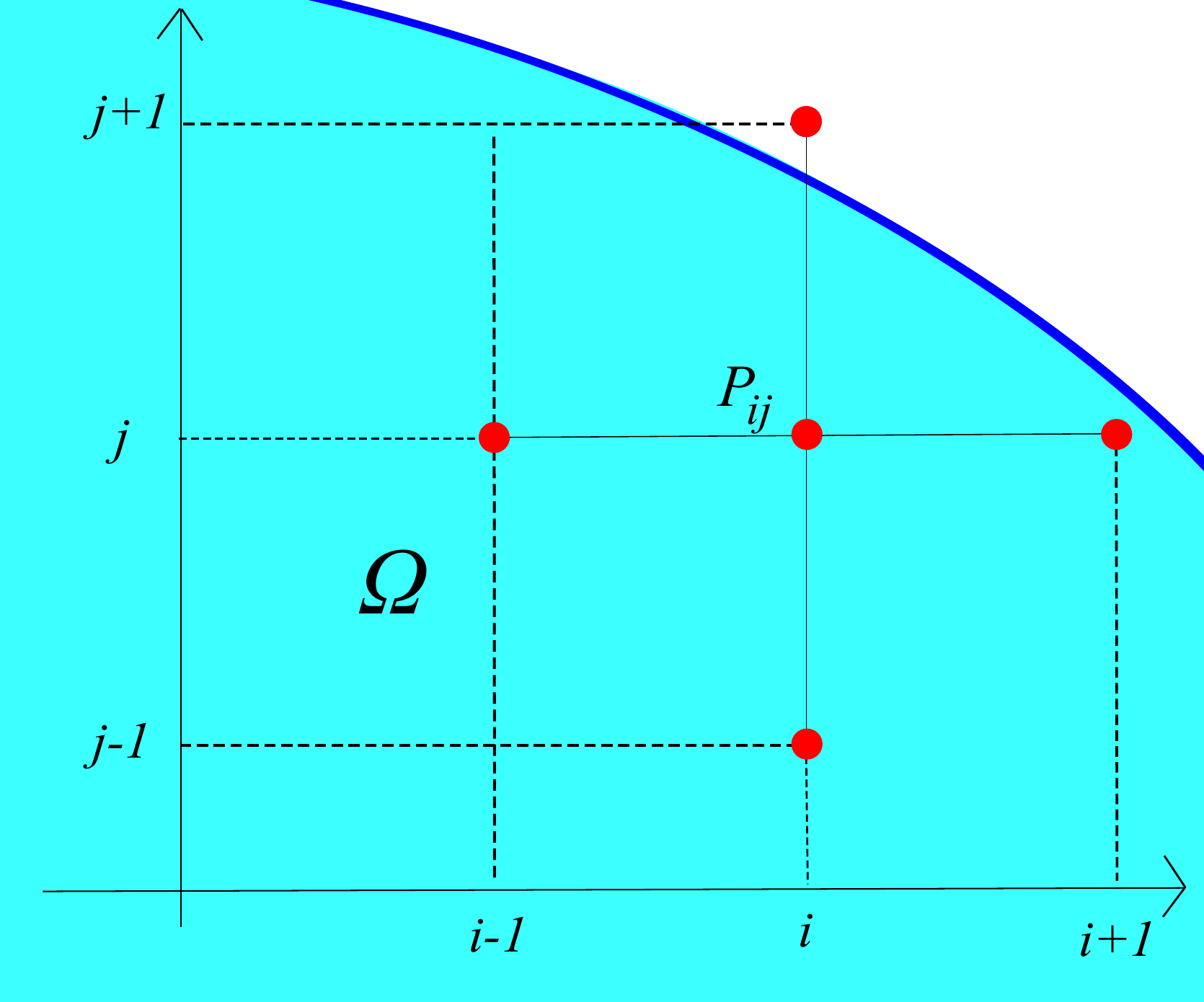}
		\caption{ \footnotesize{Five-point stencil adopted for the finite-difference discretization on inside grid points \eqref{FDs}. } }
	\label{fig:5pt_stencil}
 \end{minipage}
\end{figure}

We say that a grid point $(x_i,y_j)$ is a \textit{ghost point} if and only if both of the following conditions are satisfied:
\begin{equation}\label{ghostdef}
(x_i,y_j) \notin \Omega_h,
\qquad
\left\{ (x_i\pm h, y_j), (x_i, y_j\pm h) \right\}  \cap \Omega_h \neq \emptyset.
\end{equation}
In other words, a ghost point is a grid point that is outside the domain $\Omega$ and that has one of its four neighbour grid points inside $\Omega$. 
We call $\Gamma_h$ the set of ghost points (see Fig.~\ref{fig:domain}).


Let $\mathcal{I}_\Omega = \left\{ (i,j) \colon (x_i,y_j) \in \Omega_h \right\}$ and $\mathcal{I}_\Gamma = \left\{ (i,j) \colon (x_i,y_j) \in \Gamma_h \right\}$
be the sets of inside indices and ghost indices, respectively, and $N_i = \left| \mathcal{I}_\Omega \right|$ and $N_g = \left| \mathcal{I}_\Gamma \right|$ be their cardinalities.

We aim at approximating the solution in $\Omega_h \cup \Gamma_h$ for any time step. Therefore,
the numerical solution at each time step $t_n = n \, \Delta t$ is expressed by a vector $W^{(n)} = (\vec{s}^{(n)},\vec{c}^{(n)}) \in \R^{2(N_i+N_g)}$, whose components are $s^{(n)}_{ij}$ and $c^{(n)}_{ij}$, with $(i,j)$ varying in $\mathcal{I}_\Omega \cup \mathcal{I}_\Gamma$.

The components of $W$ are ordered by choosing a mapping
\begin{equation}\label{maprij}
\mathcal{M} \colon \left\{1,2,\ldots, N_i+N_g \right\}    \longrightarrow    \mathcal{I}_\Omega \cup \mathcal{I}_\Gamma.
\end{equation}
For the purpose of describing the numerical method, we order all $s$ variables before the $c$ ones, i.e.~$W^{(n)} = (\vec{s}^{(n)},\vec{c}^{(n)}) \in \R^{2(N_i+N_g)}$. Of course, the actual layout of the vector in the computer memory will be chosen to achieve optimal efficiency, e.g. like $W=(s_1,c_1,s_2,c_2,\ldots,s_{N_i},c_{N_i},\ldots,s_{N_i+N_g},c_{N_i+N_g})$.

In order to compute $W^{(n+1)}$ from $W^{(n)}$, a system of $2(N_i+N_g)$ nonlinear equations in $2(N_i+N_g)$ unknowns must be solved at each time step. 
The nonlinear system is obtained as follows. 
The $2 N_i$ nonlinear equations related to inside grid points are obtained by discretizing the first two equations of \eqref{discTime} on inside grid points $(x_i,y_j) \in \Omega$ using the standard five-point stencil (Fig.~\eqref{fig:5pt_stencil}) finite-difference scheme:
\begin{equation}\label{FDs}
\displaystyle \frac{ \phi \left(c_{ij}^{(n+1)} \right)s_{ij}^{(n+1)} - \phi \left(c_{ij}^{(n)} \right)s_{ij}^{(n)}}{\Delta t} -
\displaystyle \frac{\L^s_h (\vec{s}^{(n)},\vec{c}^{(n)}) + \L^s_h (\vec{s}^{(n+1)},\vec{c}^{(n+1)})}{2}
= 0,
\end{equation}
\begin{equation}\label{FDc}
\displaystyle \frac{c_{ij}^{(n+1)}-c_{ij}^{(n)}}{\Delta t} -
\displaystyle \frac{\L^c_h (\vec{s}^{(n)},\vec{c}^{(n)}) + \L^c_h (\vec{s}^{(n+1)},\vec{c}^{(n+1)})}{2} = 0
\end{equation}
where
\[
\L^s_h (\vec{s},\vec{c}) = -\frac{a}{m_c} \phi(c_{ij}) \, s_{ij} \, c_{ij}+
\frac{d}{2}\, \sum_{(i^*,j^*) \in \mathcal{N}_{ij}} 
\left( \phi \left(c_{i^*j^*}\right) + \phi \left(c_{ij} \right) \right) \left( s_{i^*j^*}  - s_{ij} \right)
\]
\[
\L^c_h (\vec{s},\vec{c}) = - \frac{a}{m_s} \phi \left(c_{ij} \right) \, s_{ij}\, c_{ij},
\]

and $\mathcal{N}_{ij}$ is the set of four neighbouring index pairs for $(i,j)$, namely:
\[
\mathcal{N}_{ij} = \left\{ (i\pm 1, j), (i, j \pm 1)  \right\}.
\]
Eq.~\eqref{FDc} is discretized also on ghost points (provided that the initial conditions for $c$ and $s$ are extrapolated outside the domain with second order accuracy, using for example the extrapolation technique proposed in~\cite{aslam2004partial}), leading to additional $N_g$ nonlinear equations. The remaining $N_g$ equations are obtained by enforcing the boundary condition for $s$ on $\Gamma$
\begin{equation}\label{bc2}
s^{(n+1)} = s_b \mbox{ on } \partial \Omega
\end{equation}
using a ghost point extrapolation technique that was already successfully adopted in other contexts (elliptic~\cite{CocoRusso:Elliptic, coco2018second, coco2012second, CCDR:LinearElasticity, coco2016hydro} and hyperbolic~\cite{chertock2018second, CocoRusso:Hyp2012} equations) and described as follows.
Let $G_{ij}=(x_i,y_j) \in \Gamma_h$ be a ghost point and $\vec{n}_{ij}=(n_x,n_y)$ be the approximated outward unit normal vector computed by a central finite-difference discretization of Eq.\ \eqref{ncurv}:
\begin{equation}\label{nij}
n_x = \frac{\tilde{n}_x}{\sqrt{\tilde{n}_x^2+\tilde{n}_y^2}}, \quad n_x = \frac{\tilde{n}_y}{\sqrt{\tilde{n}_x^2+\tilde{n}_y^2}}, \quad \mbox{  with  } \quad 
\tilde{n}_x = \frac{\varphi_{i+1,j}-\varphi_{i-1,j}}{2\, h}, \quad \tilde{n}_y = \frac{\varphi_{i,j+1}-\varphi_{i,j-1}}{2 \, h}, \quad \varphi_{ij} = \varphi(x_i,y_j).
\end{equation}
Let $\text{ST}_{ij}$ be the nine-point stencil in the upwind direction with respect to $\vec{n}_{ij}$. More precisely,
\[
\text{ST}_{ij} = \left\{ (x_i,y_j) -h\,(m_x \, k_x, m_y \, k_y), (k_x,k_y) \in \left\{ 0,1,2 \right\}^2 \right\},
\]
where $m_x=\text{SIGN} (n_x)$ and $m_y=\text{SIGN} (n_y)$, with $\text{SIGN}(x)=1$ or $-1$ if $x > 0$ or $x < 0$, respectively, and (conventionally) $\text{SIGN}(0)=1$. Since this nine-point stencil is taken in the upwind direction with respect to the normal $\vec{n}_{ij}$, it can be easily proven that $\text{ST}_9 \subseteq \Omega_h \cup \Gamma_h$, provided that the grid is sufficiently fine (i.e.\ $h$ is sufficiently small). 

\begin{figure}[!hbt]
 \begin{minipage}[c]{0.40\textwidth}
   	\centering
   	\captionsetup{width=0.80\textwidth}
		\includegraphics[width=0.79\textwidth]{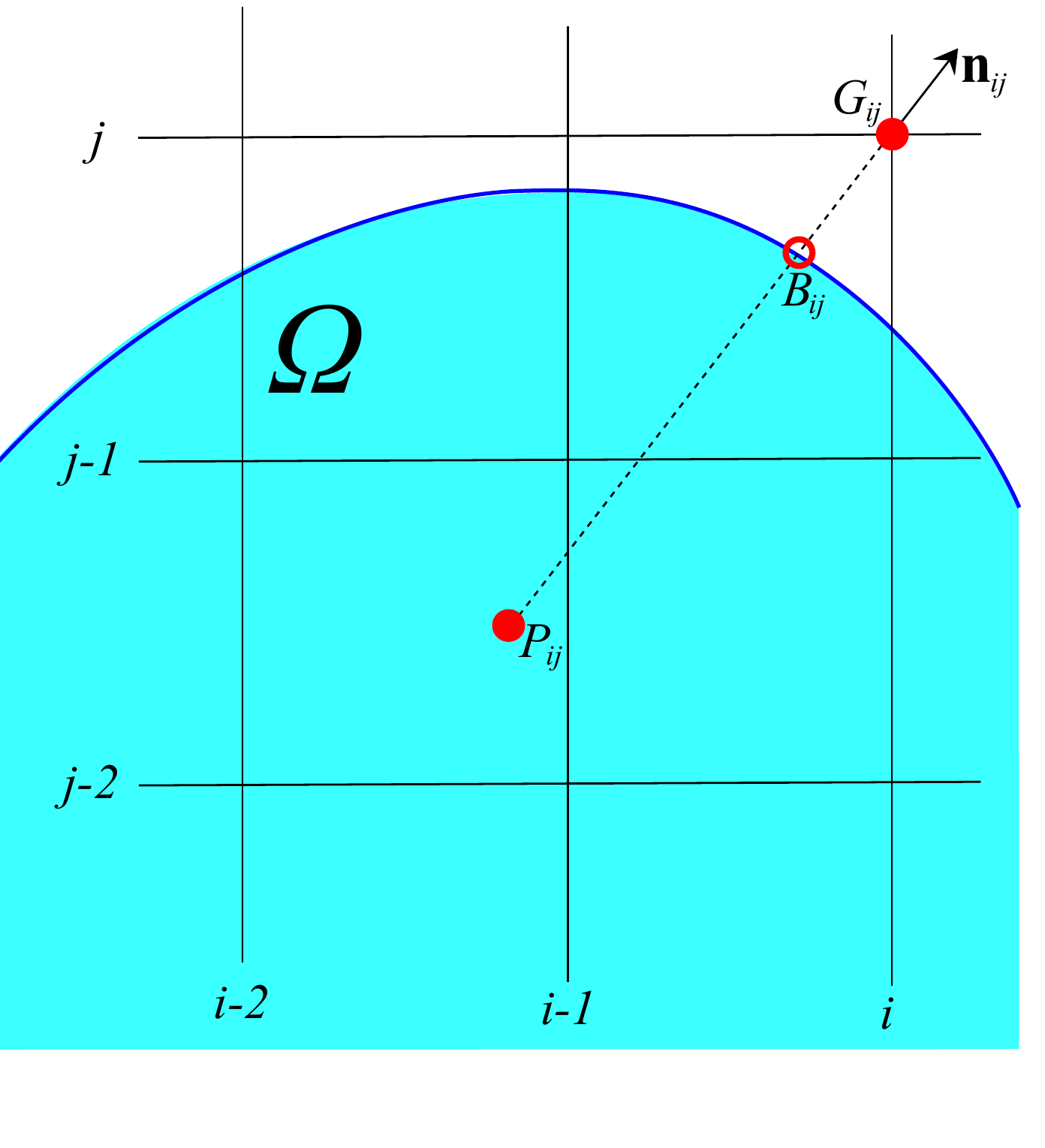}
		\caption{ \footnotesize{Boundary point $B_{ij}$ obtained from the ghost point $G_{ij}$ by solving $\varphi \left( P_{ij}+\alpha \left( G_{ij} -P_{ij} \right) \right)=0$ for $\alpha \in [0,1]$, where $P_{ij} = G_{ij} - 2 h \, \vec{n}_{ij} $.} }
	\label{fig:closestpoint}
 \end{minipage}
 \ \hspace{2mm} \hspace{3mm} \
 \begin{minipage}[c]{0.40\textwidth}
   	\centering
   	\captionsetup{width=0.80\textwidth}
		\includegraphics[width=0.99\textwidth]{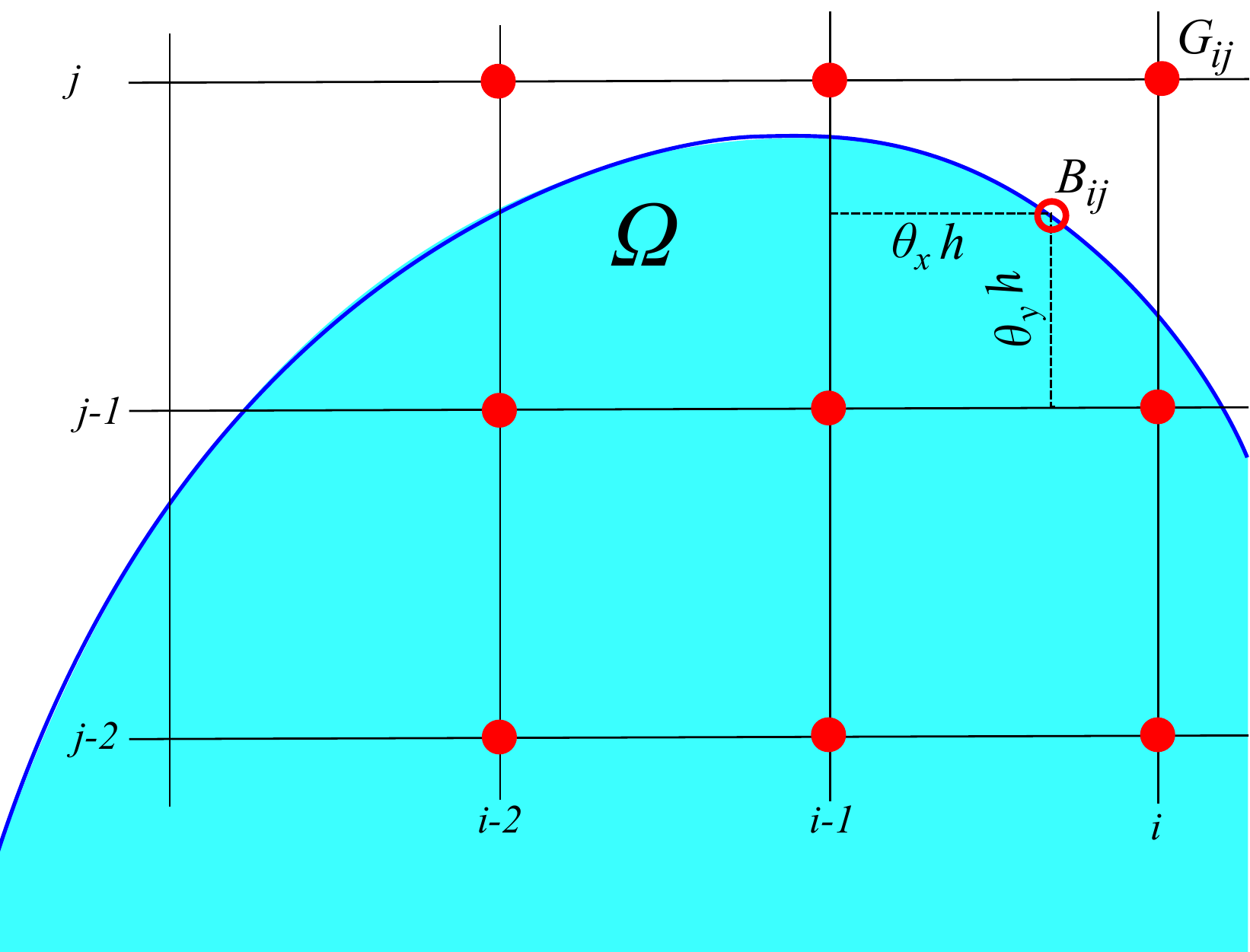}
		\caption{ \footnotesize Nine-point stencil $\text{ST}_9$ (red filled circles) associated with the ghost point $G_{ij}$ to discretize the boundary condition on $B_{ij}$ by Eq.~\eqref{JssD}. }
	\label{fig:Omega2d3}
 \end{minipage}
\end{figure}

The linear equations are obtained by prescribing the boundary condition~\eqref{bc2} on 
$B_{ij} \in \Gamma$, where $B_{ij}$ is the normal projection of $G_{ij}$ onto $\Gamma$, obtained by the following algorithm (see Fig.~\ref{fig:closestpoint}):
\begin{itemize}
\item compute $P_{ij} = G_{ij} - 2 h \, \vec{n}_{ij} $ (we have $P_{ij} \in \Omega$ provided that the grid is sufficiently fine, i.e.\ $h$ is sufficiently small);
\item apply the bisection method to solve $\tilde{\varphi} \left( P_{ij}+\alpha \left( G_{ij} -P_{ij} \right) \right)=0$ in the unknown $\alpha \in [0,1]$ (with the tolerance criterion: 
$\min \left\{ |\tilde{\varphi} \left( P_{ij}+\alpha_k \left( G_{ij} -P_{ij} \right) \right)  |, |\alpha_{k} - \alpha_{k-1}| \right\} < 10^{-6}$), where $\tilde{\varphi}$ is a biquadratic interpolation of $\varphi$ on the stencil $\text{ST}_{ij}$;
\item compute $B_{ij}=P_{ij}+\alpha \left( G_{ij} -P_{ij} \right)$.
\end{itemize}

Finally, the linear equations are:
\begin{equation}\label{BCs}
s_b-\tilde{s}^{(n+1)}(B_{ij})=0,
\end{equation}
where $\tilde{s}^{(n+1)}$
is the biquadratic interpolations of $\vec{s}^{(n+1)}$
on the stencil $\text{ST}_{ij}$.

We observe that when only Dirichlet boundary conditions are considered, other boundary discretization approaches can result more efficient and easier to implement, such as Shortley-Weller methods~\cite{Shortley-Weller:discretization} or the symmetric positive-definite discretization proposed in~\cite{Gibou:Ghost}. However, the method proposed here is robust for all types of boundary conditions (Dirichlet, Neumann, Robin, mixed) and an extension to Neumann boundary conditions is proposed in \S\ref{ssec:numtest:neumann} and tested in \S\ref{ssec:numtest:accuracy} ({\sc Test 1N and Test 2N}), making the approach suitable for an extension to more complex models of monument conservation, such as the free boundary models proposed in 1D in~\cite{CFN:08:swelling}, where a non-homogenous free boundary conditions involving fluxes is proposed to account for the swelling of the external gypsum layer, The 3D implementation is currently under investigation.

\subsection{Newton-Raphson method}\label{ssec:method:newton}
Eqs.~\eqref{FDs}, \eqref{FDc}, \eqref{BCs} 
constitute the system of nonlinear equations that we need to solve in order to advance in time from $t_n$ to $t_{n+1}$. Let $F$ be an operator $F \colon \R^{2(N_i+N_g)} \longrightarrow \R^{2(N_i+N_g)}$ such that the nonlinear system (\eqref{FDs}, \eqref{FDc}, \eqref{BCs}) can be represented by $F(W^{(n+1)})=0$. This system is solved by the Newton-Raphson method, i.e.\ by the following iterative scheme: 
\begin{enumerate}
\item take the initial guess $W^{(n+1,0)}=W^{(n)}$;
\item for $k=1,\ldots$, repeat the steps:
\begin{enumerate}
\item solve the linear system
\begin{equation}\label{NewtonLS}
J_F(W^{(n+1,k)}) \cdot \Delta W = F(W^{(n+1,k)}),
\end{equation}
where $J_F(W^{(n+1,k)})$ is the Jacobian matrix of $F$;
\item update the current guess $W^{(n+1,k+1)} = W^{(n+1,k)} - \Delta W$.
\end{enumerate}
until a suitable tolerance is reached, e.g.:
\[
\min \left\{ \left\| F(W^{(n+1,k+1)}) \right\|_\infty, \displaystyle \frac{ \left\| W^{(n+1,k+1)}-W^{(n+1,k)} \right\|_\infty }{\left\| W^{(n+1,k)} \right\|_\infty}  \right\} < 10^{-9}.
\]
\end{enumerate}
The Jacobian matrix can be represented in the compact form
\begin{equation}\label{jacob4}
J_F(W^{(n+1,k)}) = 
\begin{pmatrix}
J^{ss}(W^{(n+1,k)}) & J^{sc}(W^{(n+1,k)}) \\
J^{cs}(W^{(n+1,k)}) & J^{cc} (W^{(n+1,k)})
\end{pmatrix},
\end{equation}
where $\left\{ J^{ss},J^{sc},J^{cs},J^{cc} \right\} (W^{(n+1,k)})$ are the four $(N_i+N_g) \times (N_i+N_g)$ matrices detailed below.
Note that \eqref{jacob4} describes the logical structure of the Jacobian matrix, and that its actual layout in the computer memory may be different and will match the layout chosen for the vectors $W^{(n)}$.

We represent the rows of a $(N_i+N_g) \times (N_i+N_g)$ matrix by $3 \times 3$ stencils in the following way:
we say that the stencil
 \[
 \begin{bmatrix}
 a_{-1,1}  & a_{0,1} &  a_{1,1} \\
  a_{-1,0}  & \boxed{a_{0,0}} &  a_{1,0} \\
   a_{-1,-1}  & a_{0,-1} &  a_{1,-1}
 \end{bmatrix}
 \]
 represents the $r$-th row of a $(N_i+N_g) \times (N_i+N_g)$ matrix $A$, with $\mathcal{M}(r)=(i,j)$, if for any $1 \leq q \leq N_i + N_g$ we have that
 \[
 A_{rq} = \left\{
 \begin{matrix}
 a_{k_i,k_j} & \mbox{ if } \mathcal{M}(q) = (i+k_i,j+k_j), \mbox{ with } k_i, k_j \in \left\{ -1,0,1 \right\} \\
 0 & \mbox{ otherwise. }
  \end{matrix}
  \right.
 \]
 We have boxed the central element $a_{0,0}$ to emphasize that this element is on the main diagonal of the matrix $A$, i.e.~$A_{rr}=a_{0,0}$. Sometimes (in particular for ghost points) it may be easier to have the element of the main diagonal not necessarily at the center of the $3 \times 3$ stencil. Therefore, we use for instance the stencil
 \[
 \begin{bmatrix}
 \boxed{a_{0,0}}  & a_{1,0} &  a_{2,0} \\
 a_{0,-1}  & a_{1,-1} &  a_{2,-1} \\
 a_{0,-2}  & a_{1,-2} &  a_{2,-2} 
 \end{bmatrix}
 \]
if we want to represent the $r$-th row of $A$ as:
  \[
 A_{rq} = \left\{
 \begin{matrix}
 a_{k_i,k_j} & \mbox{ if } \mathcal{M}(q) = (i+k_i,j+k_j), \mbox{ with } k_i \in \left\{ 0,1,2 \right\}, k_j \in \left\{ -2,-1,0 \right\} \\
 0 & \mbox{ otherwise. }
  \end{matrix}
  \right.
 \]
 
\begin{figure}[!hbt]
   	\centering
   	\captionsetup{width=0.80\textwidth}
		\includegraphics[width=0.29\textwidth]{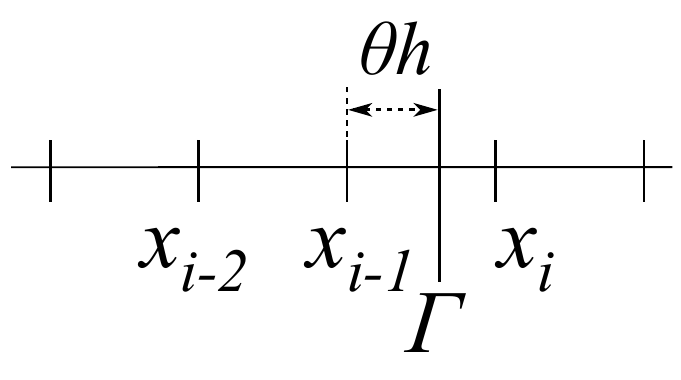}
		\caption{ \footnotesize{One-dimensional interpolation on $\Gamma$ using nodes $x_{i-2}$, $x_{i-1}$ and $x_{i}$ (Eqs.~\eqref{intcoeffsD} and \eqref{intcoeffsN}).} }
	\label{fig:Omega1d}
\end{figure}

\paragraph{Representation of the matrix $J^{ss}(W^{(n+1,k)})$.} Using this notation, if $\mathcal{M}(r) \in \mathcal{I}_\Omega$ the $r$-th row of
$J^{ss}(W^{(n+1,k)})$ is:
\[
 \begin{bmatrix}
0  & 0 &  0 \\
\\
0 &
\boxed{\frac{\phi \left(c_{ij}^{(n+1)} \right)}{\Delta t} + \displaystyle  \frac{a}{2\,m_c} \phi \left(c_{ij}^{(n+1)} \right) \, c_{ij}^{(n+1)}}
   & 
0
    \\
    \\
0  & 0 &  0
 \end{bmatrix}
\]
\begin{equation}\label{JssI}
 + \frac{d}{4}
 \begin{bmatrix}
0  & -\displaystyle  \phi \left(c_{i,j+1}^{(n+1)}\right) + \phi \left(c_{ij}^{(n+1)} \right) &  0 \\
\\
-\displaystyle \phi \left(c_{i-1,j}^{(n+1)}\right) + \phi \left(c_{ij}^{(n+1)} \right) &
\boxed{ \sum_{(i^*,j^*) \in \mathcal{N}_{ij}} 
\left( \phi \left(c_{i^*j^*}^{(n+1)}\right) + \phi \left(c_{ij}^{(n+1)} \right) \right)}
   & 
   -\displaystyle \phi \left(c_{i+1,j}^{(n+1)}\right) + \phi \left(c_{ij}^{(n+1)} \right)
    \\
    \\
0  & -\displaystyle \frac{d}{2} \left( \phi \left(c_{i,j-1}^{(n+1)}\right) + \phi \left(c_{ij}^{(n+1)} \right) \right) &  0
 \end{bmatrix}.
\end{equation}
Let $\mathcal{M}(r) \in \mathcal{I}_\Gamma$ and $n_x\geq0$, $n_y \geq 0$ (see Fig.~\ref{fig:Omega2d3} and Eq.~\eqref{nij}). Referring to
Fig.~\ref{fig:Omega1d}, observe that the three coefficients 
\begin{equation}\label{intcoeffsD}
 \begin{bmatrix}
 \displaystyle \frac{\theta (\theta-1)}{2} , \quad  (1-\theta) (1+\theta),  \quad   \frac{\theta (1+\theta)}{2} 
 \end{bmatrix}
\end{equation}
are the 1D quadratic interpolation coefficients on $\Gamma$ for grid points $x_{i-2}$, $x_{i-1}$, $x_i$, respectively, with $\theta=(\Gamma-x_{i-1})/h$.

The 2D biquadratic interpolation \eqref{JssD} is obtained as dimension by dimension quadratic interpolations. Therefore, the $r$-th row of
$J^{ss}(W^{(n+1,k)})$ is represented by
\[
 \begin{bmatrix}
\boxed{ \frac{\theta_y (1+\theta_y)}{2}}  \\
\\
 (1-\theta_y) (1+\theta_y) \\
 \\
\displaystyle    \frac{\theta_y (\theta_y-1)}{2}
 \end{bmatrix}
 \cdot
 \begin{bmatrix}
    \displaystyle \frac{\theta_x (\theta_x-1)}{2} &  (1-\theta_x) (1+\theta_x) & \boxed{ \frac{\theta_x (1+\theta_x)}{2}}
 \end{bmatrix}
 = 
 \]
 \begin{equation}\label{JssD}
 \begin{bmatrix}
    \displaystyle \frac{\theta_x (\theta_x-1) \theta_y (1+\theta_y)}{4} & \displaystyle \frac{(1-\theta_x) (1+\theta_x) \theta_y (1+\theta_y)}{2} &
  \boxed{ \frac{\theta_x (1+\theta_x) \theta_y (1+\theta_y)}{4} }\\
\\
 \displaystyle \frac{\theta_x (\theta_x-1) (1-\theta_y) (1+\theta_y)}{2}  & (1-\theta_x) (1+\theta_x) (1-\theta_y) (1+\theta_y) &
\displaystyle \frac{\theta_x (1+\theta_x) (1-\theta_y) (1+\theta_y)}{2}  \\
 \\
   \displaystyle \frac{\theta_x (\theta_x-1) \theta_y (\theta_y-1)}{4} & \displaystyle \frac{(1-\theta_x) (1+\theta_x) \theta_y (\theta_y-1)}{2} &
  \displaystyle \frac{\theta_x (1+\theta_x) \theta_y (\theta_y-1)}{4} 
 \end{bmatrix},
\end{equation}
with
\[
(\theta_x, \theta_y) = \frac{|B_{ij}-(x_{i+1},y_{j-1})|}{h}.
\]
The other three possible cases  $\left\{ n_x < 0, n_y \geq 0\right\}$, $\left\{ n_x \geq 0, n_y < 0\right\}$ and $\left\{ n_x < 0, n_y < 0\right\}$ are obtained similarly. Observe that the coefficients of the other three cases are the same as in \eqref{JssD}, but in a different order, and that the boxed coefficient (i.e.\ the coefficient that will populate the main diagonal of the matrix) has the same expression in all cases.

\paragraph{Representation of the matrix $J^{sc}(W^{(n+1,k)})$.}
The $r$-th row of $J^{sc}(W^{(n+1,k)})$ is
\[
 \begin{bmatrix}
0  & 0 &  0 \\
\\
0 &
\boxed{
\displaystyle \frac{\phi' \left(c_{ij}^{(n+1)}\right) s_{ij}^{(n+1)}}{\Delta t} + 
\frac{a}{2\,m_c} \left( \phi' \left(c_{ij}^{(n+1)}\right)  \, s_{ij}^{(n+1)} \, c_{ij}^{(n+1)} + \phi \left(c_{ij}^{(n+1)}\right)  \, s_{ij}^{(n+1)}   \right) 
}
   & 
   0
    \\
    \\
0  & 0 &  0
 \end{bmatrix}
\]
\begin{equation}\label{JscI}
 + \frac{d}{4} \sum_{(i^*,j^*) \in \mathcal{N}_{ij}} \left( s_{ij}^{(n+1)}  - s_{i^*j^*}^{(n+1)} \right)
 \begin{bmatrix}
0  & \phi' \left(c_{i,j+1}^{(n+1)} \right)  &  0 \\
\\
\phi' \left(c_{i-1,j}^{(n+1)} \right) &
\boxed{
\phi' \left(c_{i,j}^{(n+1)} \right) 
}
   & 
\phi' \left(c_{i+1,j}^{(n+1)} \right)
    \\
    \\
0  & \phi' \left(c_{i,j-1}^{(n+1)} \right) &  0
 \end{bmatrix}
\end{equation}
if $\mathcal{M}(r) \in \mathcal{I}_\Omega$, or a null row if $\mathcal{M}(r) \in \mathcal{I}_\Gamma$.

\paragraph{Representation of the matrix $J^{cs}(W^{(n+1,k)})$.}
The $r$-th row of $J^{cs}(W^{(n+1,k)})$ is
\begin{equation}\label{JcsI}
 \begin{bmatrix}
0  & 0 &  0 \\
\\
0 &
\boxed{
\frac{a}{m_s} \phi \left(c_{ij}^{(n+1)} \right) \, s_{ij}^{(n+1)}\, c_{ij}^{(n+1)}
}
   & 
   0
    \\
    \\
0  & 0 &  0
 \end{bmatrix}
\end{equation}
both if $\mathcal{M}(r) \in \mathcal{I}_\Omega$ or $\mathcal{M}(r) \in \mathcal{I}_\Gamma$.

\paragraph{Representation of the matrix $J^{cc}(W^{(n+1,k)})$.}
The $r$-th row of $J^{cc}(W^{(n+1,k)})$ is represented by
\begin{equation}\label{JccI}
 \begin{bmatrix}
0  & 0 &  0 \\
\\
0 &
\boxed{
\displaystyle \frac{1}{\Delta t} +  
\frac{a}{m_s} \left( \phi' \left(c_{ij}^{(n+1)}\right)  \, s_{ij}^{(n+1)} \, c_{ij}^{(n+1)} + \phi \left(c_{ij}^{(n+1)}\right)  \, s_{ij}^{(n+1)}   \right) 
}
   & 
   0
    \\
    \\
0  & 0 &  0
 \end{bmatrix}.
\end{equation}
both if $\mathcal{M}(r) \in \mathcal{I}_\Omega$ or $\mathcal{M}(r) \in \mathcal{I}_\Gamma$.

Summarizing, the four sub-matrices of the Jacobian matrix \eqref{jacob4} can be represented in matrix form (see Eqs.~\eqref{JssI}, \eqref{JssD}, \eqref{JscI}, \eqref{JcsI}, \eqref{JccI}):
\[
J^{ss}(W^{(n+1,k)}) = D_{ss}(W^{(n+1,k)}) + \frac{d}{4} M_{ss}(W^{(n+1,k)}) + R_{ss},
\quad 
J^{sc}(W^{(n+1,k)}) = D_{sc}(W^{(n+1,k)}) + \frac{d}{4} M_{sc}(W^{(n+1,k)}),
\]
\[
J^{cs}(W^{(n+1,k)}) = D_{cs}(W^{(n+1,k)}),
\quad 
J^{cc}(W^{(n+1,k)}) = D_{cc}(W^{(n+1,k)}), 
\]
where $D_{ss}(W^{(n+1,k)})$, $D_{sc}(W^{(n+1,k)})$, $D_{cs}(W^{(n+1,k)})$, $D_{cc}(W^{(n+1,k)})$ are four $(N_i+N_g) \times (N_i+N_g)$ diagonal matrices, $M_{ss}(W^{(n+1,k)})$, $M_{sc}(W^{(n+1,k)})$ are two $(N_i+N_g) \times (N_i+N_g)$ penta-diagonal matrices, while $R_{ss}$ is a $(N_i+N_g) \times (N_i+N_g)$ nine-diagonal matrix that does not depend on $W^{(n+1,k)}$ and then can be precomputed at the beginning of the numerical simulation.

\subsection{Multigrid method}\label{ssec:method:multigrid}
The linear system \eqref{NewtonLS} is solved by a multigrid tecnique, as described in this section. In particular, we will introduce the relaxation operator (\S\ref{relop}) and the transfer (restriction and interpolation) operators (\S\ref{sec:restrictionOp} and \ref{sec:interpolationOp}). The multigrid method can then be easily implemented from these operators (we refer the reader to any book on multigrid methods for a comprehensive presentation, such as~\cite{Trottemberg:MG}).
In this paper we implement the $W-$cycle multigrid scheme and compare its convergence factor against the one predicted by the Local Fourier Analysis for $W-$cycle multigrid schemes in rectangular domains. Extensions to multigrid schemes more efficient than $W-$cycle such as Full-multigrid~\cite[Ch.~2.6]{Trottemberg:MG} can be also easily implemented following the same approach proposed in this paper.

\subsubsection{Relaxation operator}\label{relop}
Let us write the linear system \eqref{NewtonLS} as:
\begin{equation*}
\begin{bmatrix}
J^{ss}(W^{(n+1,k)}) & J^{sc}(W^{(n+1,k)}) \\
J^{cs}(W^{(n+1,k)}) & J^{cc}(W^{(n+1,k)}) 
\end{bmatrix}
\cdot 
\begin{bmatrix}
\Delta \vec{s} \\
\Delta \vec{c}
\end{bmatrix}
= 
\begin{bmatrix}
F^s(W^{(n+1,k)}) \\
F^c(W^{(n+1,k)})
\end{bmatrix}.
\end{equation*}
In order to have an efficient multigrid method,
the relaxation operator must satisfy the smoothing property, i.e.~the high frequency components of the defect should be dumped quickly after few relaxations, in such a way that the residual linear system is well represented on a coarser grid (see ~\cite[Ch.~2.1]{Trottemberg:MG}). Well known relaxation operators that satisfy the smoothing property for scalar elliptic equations are Gauss-Seidel and weighted-Jacobi (with $\omega=4/5$). It is known that the classical Gauss-Seidel scheme may underperform for systems of partial differential equations, where the collective Gauss-Seidel scheme is usually preferred (\cite[Ch.~8]{Trottemberg:MG}).
In this paper we implement a collective Gauss-Seidel scheme, i.e.~a $2\times 2$ linear system is solved at each internal grid point to update $\Delta s_{ij}$ and $\Delta c_{ij}$ simultaneously, while an appropriate relaxation is performed on ghost points (\cite{CocoRusso:Elliptic, coco2018second, coco2012second}). 
In detail, the relaxation scheme spans all $r = 1,\ldots, N_i+N_g$. If $\mathcal{M} (r) = (i,j) \in \mathcal{I}_\Omega$, then  $\Delta s_{ij}$ and $\Delta c_{ij}$ are updated as follows:
\begin{eqnarray*}
\begin{bmatrix}
\left( \Delta s^{(m+1)} \right)_r \\
\left( \Delta c^{(m+1)} \right)_r 
\end{bmatrix}
=
\begin{bmatrix}
\left( \Delta s^{(m)} \right)_r \\
\left( \Delta c^{(m)} \right)_r 
\end{bmatrix}
+ 
P^{-1}
\cdot 
\begin{bmatrix}
F^s (W^{(n+1,k)})_{r}- \left( J^{ss} (W^{(n+1,k)})_r \cdot \Delta \vec{s}^{*} + J^{sc}(W^{(n+1,k)})_r \cdot \Delta \vec{c}^{*} \right) \\
F^c (W^{(n+1,k)})_{r}- \left( J^{cs} (W^{(n+1,k)})_r \cdot \Delta \vec{s}^{*} + J^{cc}(W^{(n+1,k)})_r \cdot \Delta \vec{c}^{*} \right)
\end{bmatrix},
\end{eqnarray*}
where
$F^s (W^{(n+1,k)})_{r}$ and $F^c (W^{(n+1,k)})_{r}$ are the $r-$th components of the vectors $F^s (W^{(n+1,k)})$ and $F^c (W^{(n+1,k)})$, respectively, and $J^{\left\{ss,sc,cs,cc \right\}} (W^{(n+1,k)})_r$ is the $r-$th row of the matrix $J^{\left\{ss,sc,cs,cc \right\}} (W^{(n+1,k)})$ (see Eqs.~\eqref{JssI}, \eqref{JscI}, \eqref{JcsI}, \eqref{JccI}). 
We have denoted by $\Delta \vec{s}^{*} $ and $\Delta \vec{c}^{*} $ the current approximations of the Gauss-seidel iteration, i.e.~they are vectors whose $q-$th component is defined by:
\begin{equation}\label{deltasc}
\left( \Delta \vec{s}^{*} \right)_q =
\left\{
\begin{matrix}
 \left( \Delta \vec{s}^{(m+1)} \right)_q & \mbox{ if } 1\leq q < r \\
  \left( \Delta \vec{s}^{(m)} \right)_q & \mbox{ if } r \leq q \leq N_i + N_g
 \end{matrix}
 \right., \quad
 \left( \Delta \vec{c}^{*} \right)_q =
\left\{
\begin{matrix}
 \left( \Delta \vec{c}^{(m+1)} \right)_q & \mbox{ if } 1\leq q < r \\
  \left( \Delta \vec{c}^{(m)} \right)_q & \mbox{ if } r \leq q \leq N_i + N_g
 \end{matrix}
 \right. ,
\end{equation}
and 
\begin{equation}\label{matrixP}
P = \begin{bmatrix}
\left( J^{ss} (W^{(n+1,k)}) \right)_{rr} & \left( J^{sc} (W^{(n+1,k)}) \right)_{rr} \\
\left( J^{cs} (W^{(n+1,k)}) \right)_{rr} & \left( J^{cc} (W^{(n+1,k)}) \right)_{rr}
\end{bmatrix}.
\end{equation}
We observe that the classical (pointwise) Gauss-Seidel scheme can be obtained by replacing the matrix \eqref{matrixP} with
\begin{equation}\label{matrixPclassical}
P = \begin{bmatrix}
\left( J^{ss} (W^{(n+1,k)}) \right)_{rr} & 0 \\
0 & \left( J^{cc} (W^{(n+1,k)}) \right)_{rr}
\end{bmatrix}.
\end{equation}
If $\mathcal{M} (r) = (i,j) \in \mathcal{I}_\Gamma$, then $\Delta s_{ij}$ and $\Delta c_{ij}$ are updated as follows:
\begin{equation}\label{RelGhost}
\left( \Delta s^{(m+1)} \right)_{r} = \left( \Delta s^{(m)} \right)_{r} + \tau^s \left(
F^s (W^{(n+1,k)})_{r}- J^{ss} (W^{(n+1,k)})_r \cdot \Delta \vec{s}^{*} \right) \\
\end{equation}
\begin{multline*}
\left( \Delta c^{(m+1)} \right)_{r} = \left( \Delta c^{(m)} \right)_{r} \\
+ \left( J^{cc} (W^{(n+1,k)}) \right)_{rr}^{-1} \left(
F^c (W^{(n+1,k)})_{r}- \left( J^{cs} (W^{(n+1,k)})_r \cdot \Delta \vec{s}^{*} + J^{cc}(W^{(n+1,k)})_r \cdot \Delta \vec{c}^{*} \right) \right),
\end{multline*}
where $J^{\left\{ss \right\}} (W^{(n+1,k)})_r$ is the $r-$th row of the matrix $J^{\left\{ss \right\}} (W^{(n+1,k)})$ (see Eqs.~\eqref{JssD}). 
If we choose $\tau^s$ as in \eqref{matrixPclassical} to have a Gauss-Seidel iteration, i.e.~
\[
\tau^s= \left( \frac{\theta_x (1+\theta_x) \theta_y (1+\theta_y)}{4} \right)^{-1},
\]
the relaxation scheme may not converge (we observe for example that rows \eqref{JssD} are not diagonally dominant). Following the idea proposed in~\cite{CocoRusso:Elliptic, coco2018second, coco2012second}, the parameter $\tau^s$ is chosen in such a way that a proper CFL condition is satisfied for the iterations \eqref{RelGhost}. In particular, we want to ensure that the absolute value of the coefficient of $\left( \Delta (s)^{(m)} \right)_{r}$ in the right-hand side of \eqref{RelGhost} is smaller than one, i.e.~
\begin{equation}\label{CFLcond}
\begin{matrix}
\left| 1 - \tau^s   \displaystyle \frac{\theta_x (1+\theta_x) \theta_y (1+\theta_y)}{4}  \right| < 1 \\
\end{matrix}
\end{equation}
We choose a parameter $\tau^s$ that satisfies the condition \eqref{CFLcond} for any $0\leq \theta_x, \theta_y \leq 1$.
This is achieved by:
$
0 < \tau^s < 1.
$
For practical purposes, we use
$
\tau^s =0.9.
$

Finally, although in this paper we will use the lexicographic order of the map \eqref{maprij}, more efficient (collective) Gauss-Seidel schemes for multigrid methods, such as Red-Black Gauss-Seidel, can be easily implemented.

\subsubsection{Restriction operator}\label{sec:restrictionOp}
After $\nu_1$ pre-relaxation iterations (\S\ref{relop}), we compute the defects $\vec{r}^s_h$ and $\vec{r}^c_h$:
\begin{equation*}
\begin{matrix}
\vec{r}^s_h = F^s (W^{(n+1,k)})- \left( J^{ss} (W^{(n+1,k)}) \cdot \Delta \vec{s}^{(\nu_1)} + J^{sc}(W^{(n+1,k)}) \cdot \Delta \vec{c}^{(\nu_1)} \right) \\
\vec{r}^c_h = F^c (W^{(n+1,k)})- \left( J^{cs} (W^{(n+1,k)}) \cdot \Delta \vec{s}^{(\nu_1)} + J^{cc}(W^{(n+1,k)}) \cdot \Delta \vec{c}^{(\nu_1)} \right)
\end{matrix}
\end{equation*}

that will be restricted to the coarser grid (with spatial step $2h$) by a suitable restriction operator $\mathcal{I}_{2h}^h$: 
\[
\vec{r}^s_{2h} = \mathcal{I}_{2h}^{h\, (s)} \, \vec{r}^s_h, \quad \vec{r}^c_{2h} = \mathcal{I}_{2h}^{h\, (c)} \, \vec{r}^c_h.
\]
We observe that the defect $\vec{r}^s_h$ is discontinuous across the boundary, because the defect on $\mathcal{I}_\Omega$ is related to the internal equations, while its values on $\mathcal{I}_\Gamma$ are referred to the boundary conditions. For this reason, the restriction of the internal equations must use values only from $\mathcal{I}_\Omega$. (see~\cite{CocoRusso:Elliptic})
To this purpose, the restriction operator $\mathcal{I}_{2h}^{h\, (s)}$ is appropriately modified for internal grid points in the vicinity of $\Gamma = \partial \Omega$. We recall the full-weighting restriction operator (see~\cite[Ch.~2.3.3]{Trottemberg:MG}):
\begin{equation}\label{Fstencil}
\I_{2h}^h = \frac{1}{16}
\left[
\begin{array}{ccc}
1 & 2 & 1 \\
2 & 4 & 2 \\
1 & 2 & 1 \\
\end{array}
\right]_{2h}^{h}.
\end{equation}
In general, by the stencil notation
\begin{equation}\label{Fstencilgen}
\I_{2h}^h = 
\left[
\begin{array}{ccccc}
 & \vdots & \vdots & \vdots & \\
\cdots & s_{-1,-1} & s_{-1,0} & s_{-1,1} & \cdots \\
\cdots & s_{0,-1} & s_{0,0} & s_{0,1} & \cdots \\
\cdots & s_{1,-1} & s_{1,0} & s_{1,1} & \cdots \\
       & \vdots   &  \vdots &  \vdots &      
\end{array}
\right]_{2h}^{h}
\end{equation}
we will denote the restriction operator $\I_{2h}^h$ defined by:
\[
   \I_{2h}^h \wh (x,y) = \sum_{(i,j)\in R_k} s_{i,j} \wh(x+jh,y+ih),
\]
where only a finite number of coefficients $t_{i,j}$ is different from zero, and $R_k\equiv \left\{ -k,\ldots,k \right\}^2$ for some positive integer $k$. A second order restriction is achieved with $k=1$.

\begin{figure}[!hbt]%
\begin{minipage}{0.30\textwidth}
  \begin{center}
  \includegraphics[width=0.99\textwidth]{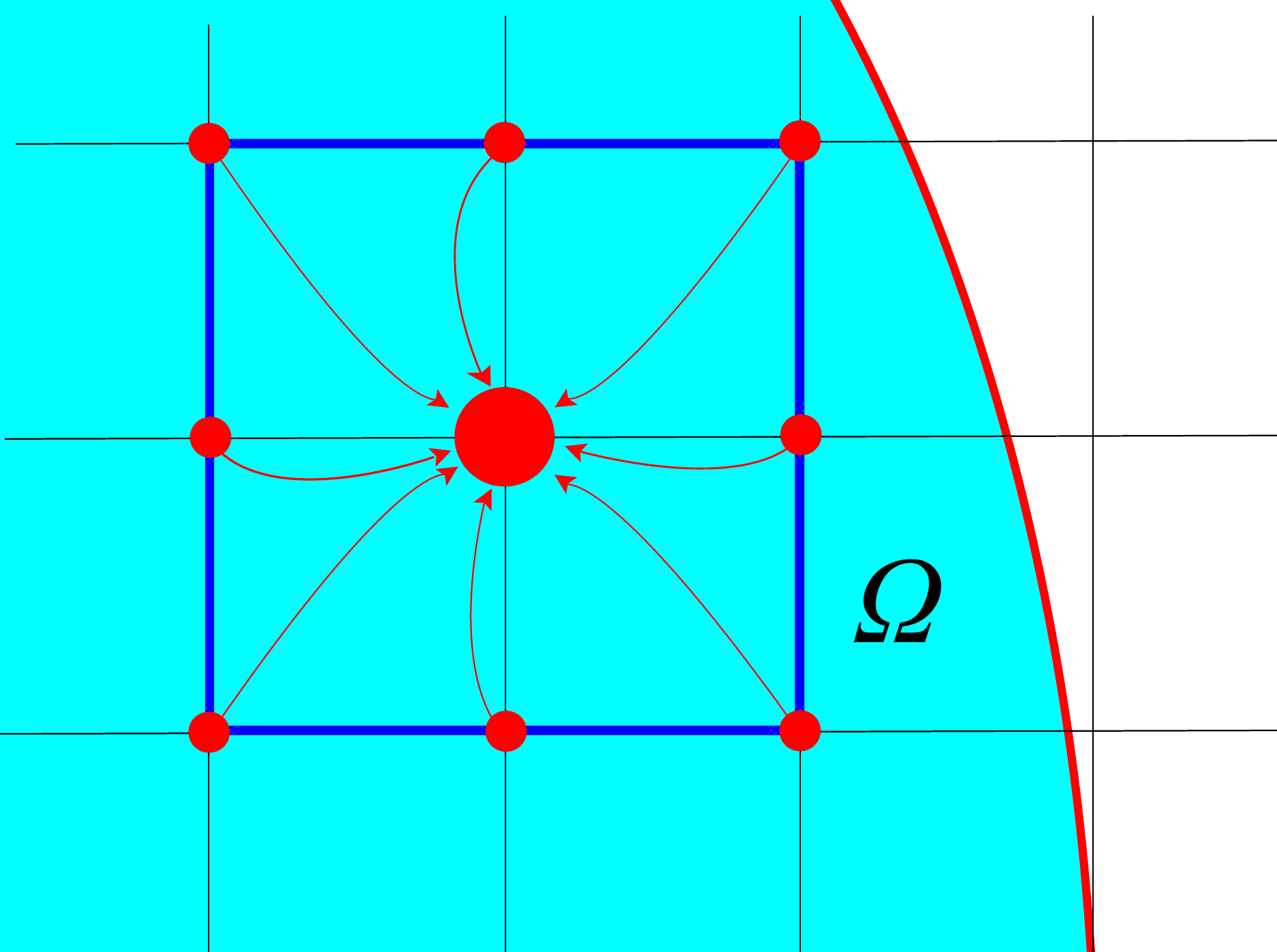}%
  \end{center}
\end{minipage}
\begin{minipage}{0.30\textwidth}
  \begin{center}
  \includegraphics[width=0.99\textwidth]{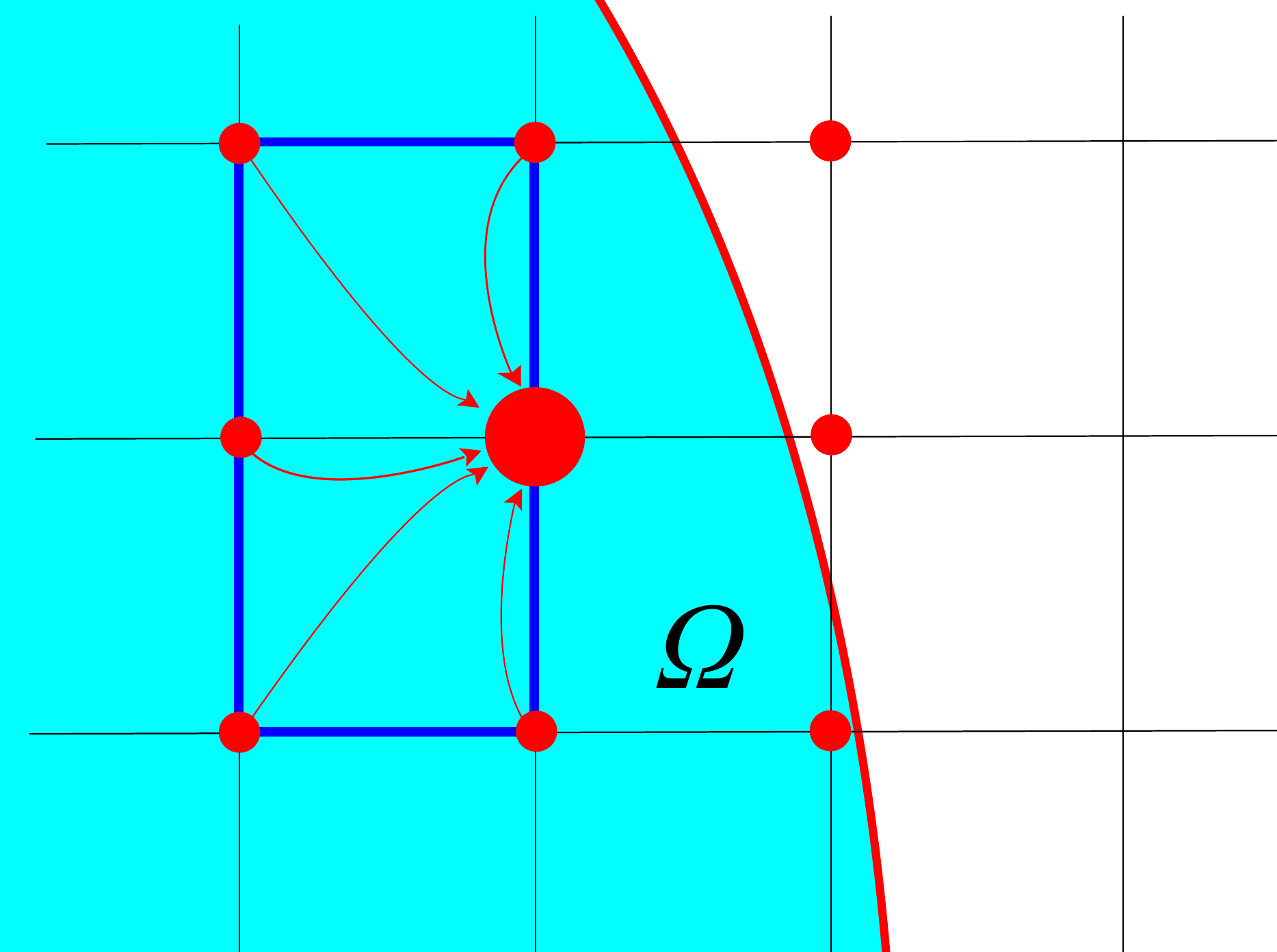}%
  \end{center}
\end{minipage}
\begin{minipage}{0.30\textwidth}
  \begin{center}
  \includegraphics[width=0.99\textwidth]{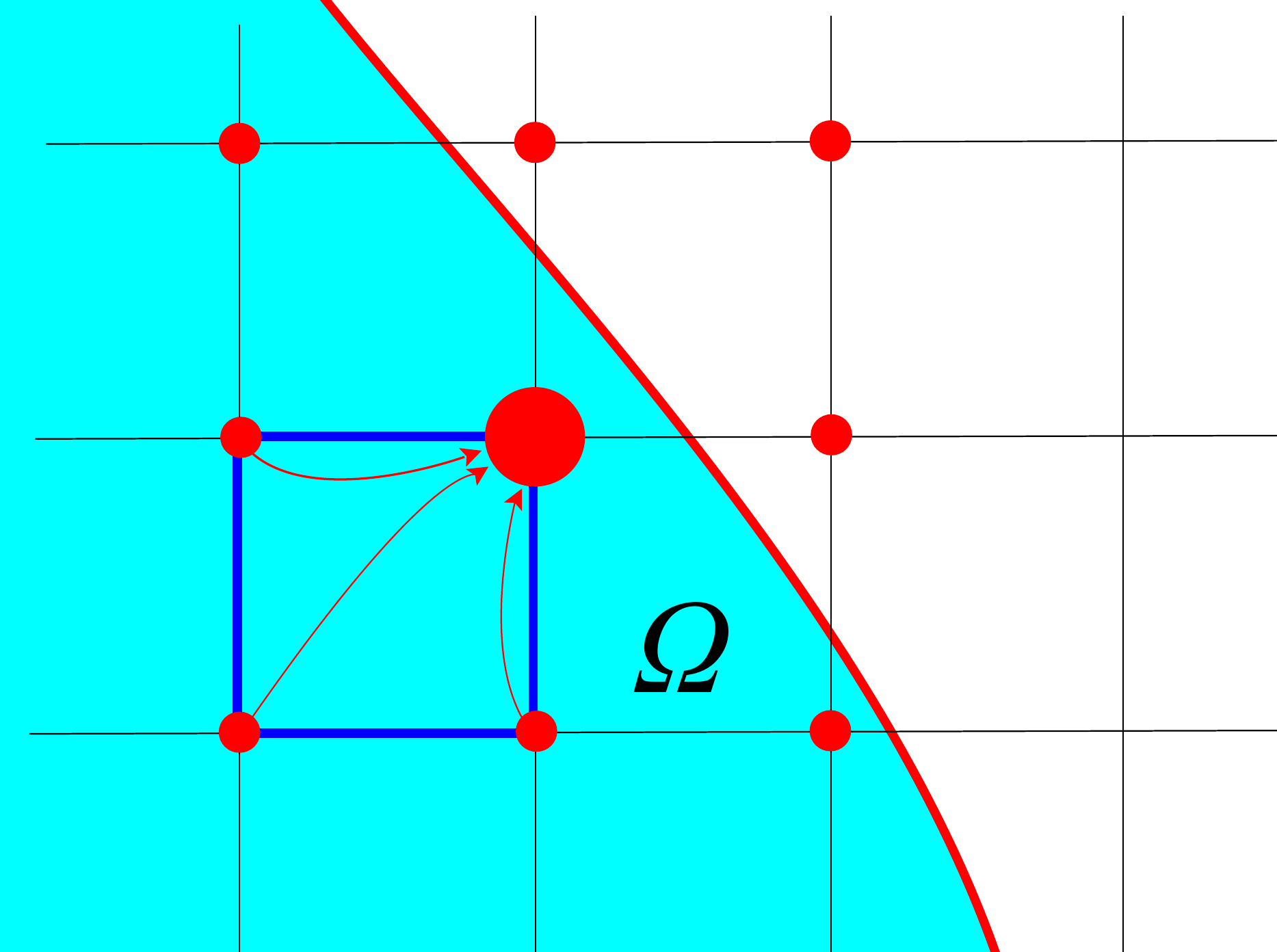}%
  \end{center}
\end{minipage}
\begin{minipage}{0.30\textwidth}
\begin{equation*}
\frac{1}{16}
\left[
\begin{array}{ccc}
1 & 2 & 1 \\
2 & 4 & 2 \\
1 & 2 & 1 \\
\end{array}
\right]_{2h}^{h}
\end{equation*}
\end{minipage}
\hspace{1mm}
\begin{minipage}{0.30\textwidth}
\begin{equation*}
\frac{1}{16}
\left[
\begin{array}{ccc}
2 & 2 & 0 \\
4 & 4 & 0 \\
2 & 2 & 0 \\
\end{array}
\right]_{2h}^{h}
\end{equation*}
\end{minipage}
\hspace{1mm}
\begin{minipage}{0.30\textwidth}
\begin{equation*}
\frac{1}{16}
\left[
\begin{array}{ccc}
0 & 0 & 0 \\
4 & 4 & 0 \\
4 & 4 & 0 \\
\end{array}
\right]_{2h}^{h}
\end{equation*}
\end{minipage}
 \captionsetup{width=0.80\textwidth}
\caption{{ Top: nine point stencil (small circles) around a grid point on the coarse grid (bold circle). The rectangle is a visual representation of the fine grid points requested in the restriction operator. Bottom: the respective stencils of the restriction operator in matrix form (Eq.~\eqref{Fstencilgen}). }}
\label{fig:Neigh}
\end{figure}

Following the same idea of~\cite{CocoRusso:Elliptic, coco2018second}, we modify the restriction operator close to the boundary in order to disregard the values of the defects on ghost points (see Fig.~\ref{fig:Neigh} and~\cite{CocoRusso:Elliptic, coco2018second} for more details).

The restriction of the boundary condition is performed using the same idea, namely using only points outside the domain (i.e.~either ghost points or inactive grid points), provided that the defect $\vec{r}^s_h$ is firstly defined in the neighbouring inactive points by extrapolating the ghost value constantly along the normal direction to the boundary $\Gamma$. This can be achieved by solving the transport equations
\begin{equation*}
\frac{\partial \vec{r}^{s}_h}{\partial \tau} + \nabla \vec{r}^{s}_h \cdot \vec{n} = 0.
\end{equation*}
for a few steps of a fictitious time $\tau$, (using for example Euler explicit method), where $\vec{n} = \nabla \varphi / |\nabla \varphi|$ is the unit normal vector.

\subsubsection{Interpolation operator}\label{sec:interpolationOp}
The defect equations
\begin{equation*}
\begin{bmatrix}
J^{ss}(W^{(n+1,k)}) & J^{sc}(W^{(n+1,k)}) \\
J^{cs}(W^{(n+1,k)}) & J^{cc}(W^{(n+1,k)}) 
\end{bmatrix}
\cdot 
\begin{bmatrix}
\vec{e}^s_{2h} \\
\vec{e}^c_{2h}
\end{bmatrix}
= 
\begin{bmatrix}
\vec{r}^s_{2h} \\
\vec{r}^c_{2h}
\end{bmatrix}
\end{equation*}
are solved recursively on the coarser grid (where $J^{\left\{ss,sc,cs,cc \right\}} (W^{(n+1,k)})$ are approximated in the coarser grid using the same technique as in the fine grid), and then the error is interpolated back to the finer grid:
\[
\vec{e}^s_{h} = \mathcal{I}^{2h}_h \vec{e}^s_{2h}, \quad \vec{e}^c_{h} = \mathcal{I}^{2h}_h \vec{e}^c_{2h}.
\]
Since errors $\vec{e}^{s}_{2h}$ and $\vec{e}^{c}_{2h}$ are continuous across the boundary, we can use values from both sides of the boundary in the same stencil and therefore we can adopt a standard linear interpolation operator for all points:
\begin{equation*}
\mathcal{I}_{h}^{2h}=
\frac{1}{4}
\left]
\begin{array}{ccc}
1 & 2 & 1 \\
2 & 4 & 2 \\
1 & 2 & 1 \\
\end{array}
\right[_{h}^{2h}.
\end{equation*}
Finally, $\nu_2$ post-relaxation iterations (\S\ref{relop}) are performed on the finer grid.

\subsection{Neumann boundary conditions}\label{ssec:numtest:neumann}
In this section we extend the numerical approach to the case of Neumann boundary conditions for $s$. The implementation of more sophisticated boundary conditions (Robin, mixed, etc.) will then be straightforward. We consider the Neumann boundary condition 
\[
\displaystyle \frac{\partial s^{(n+1)}}{\partial n} = s_b  \mbox{ on } \partial \Omega 
\]
instead of~\eqref{bc2}.

The spatial discretization~\eqref{BCs} becomes:
\begin{equation}\label{BCsN}
\left. \left( \nabla \tilde{s}^{(n+1)} \cdot  \frac{\nabla \tilde{\varphi} }{| \nabla \tilde{\varphi} |} \right) \right|_{B_{ij}} =s_b(B_{ij}),
\end{equation}
where $\tilde{s}^{(n+1)}$ and $\tilde{\varphi}$ are the the biquadratic interpolations of $\vec{s}^{(n+1)}$ and $\varphi$, respectively,
on the stencil $\text{ST}_{ij}$.
The $r$-th row of $J^{ss}(W^{(n+1,k)})$ when $\mathcal{M}(r) \in \mathcal{I}_\Gamma$ (Eq.~\eqref{JssD})
is obtained by using the condition \eqref{BCsN}:
\begin{equation}\label{Neumanntilde}
\tilde{n}_x \frac{\partial \tilde{s}}{\partial x} + \tilde{n}_y \frac{\partial \tilde{s}}{\partial y}  = 0, \quad \mbox{ with }
\tilde{n}_x = \frac{ \tilde{n}^*_x }{\sqrt{(\tilde{n}^*_x)^2+(\tilde{n}^*_y)^2}}, \quad \tilde{n}_y = \frac{\tilde{n}^*_y}{\sqrt{(\tilde{n}^*_x)^2+(\tilde{n}^*_y)^2}},
\quad \tilde{n}^*_x = \frac{\partial \tilde{\varphi}}{\partial x},
\quad \tilde{n}^*_y = \frac{\partial \tilde{\varphi}}{\partial y}.
\end{equation}
Since the coefficients of the 1D quadratic approximation of the first derivative on $\Gamma$ are (see Fig.~\ref{fig:Omega1d})
\begin{equation}\label{intcoeffsN}
\displaystyle  \frac{1}{h}
 \begin{bmatrix}
\displaystyle \theta - \frac{1}{2}, & \displaystyle -2 \theta,  &  \displaystyle \theta + \frac{1}{2} 
 \end{bmatrix}
\end{equation}
for grid points $x_{i-2}$, $x_{i-1}$, $x_i$, respectively, then the $r$-th row of $J^{ss}(W^{(n+1,k)})$ is represented by (from \eqref{Neumanntilde}, see Fig.~\ref{fig:Omega2d3}):
\[
\tilde{n}_x
 \begin{bmatrix}
\boxed{ \frac{\theta_y (1+\theta_y)}{2}}  \\
\\
 (1-\theta_y) (1+\theta_y) \\
 \\
\displaystyle    \frac{\theta_y (\theta_y-1)}{2}
 \end{bmatrix}
 \cdot 
\displaystyle  \frac{1}{h}
 \begin{bmatrix}
\displaystyle \theta_x - \frac{1}{2} & \displaystyle -2 \theta_x  &  \boxed{ \theta_x + \frac{1}{2} }  
 \end{bmatrix}
 \]
\[
+
\tilde{n}_y
\displaystyle  \frac{1}{h}
 \begin{bmatrix}
\boxed{ \theta_y + \frac{1}{2} }  \\
\\
\displaystyle -2 \theta_y \\
 \\
\displaystyle \theta_y - \frac{1}{2}
 \end{bmatrix}
 \cdot 
 \begin{bmatrix}
    \displaystyle    \frac{\theta_x (\theta_x-1)}{2} & (1-\theta_x) (1+\theta_x)  & \boxed{ \frac{\theta_x (1+\theta_x)}{2} }
 \end{bmatrix}
\]
 \[
= \frac{\tilde{n}_x }{h}
 \begin{bmatrix}
\displaystyle \frac{\theta_y (1+\theta_y)}{2} \displaystyle \left(\theta_x - \frac{1}{2}\right) &   \displaystyle  -\theta_x\theta_y (1+\theta_y) & \boxed{ \displaystyle   \frac{\theta_y (1+\theta_y)}{2} \left( \theta_x + \frac{1}{2} \right)}  \\
\\
\displaystyle (1-\theta_y) (1+\theta_y) \left(\theta_x - \frac{1}{2}\right) & \displaystyle -2 \theta_x  (1-\theta_y) (1+\theta_y) & \displaystyle (1-\theta_y) (1+\theta_y) \left( \theta_x+\frac{1}{2} \right)  \\
 \\
 \displaystyle  \frac{\theta_y (\theta_y-1)}{2} \left(\theta_x - \frac{1}{2}\right) & \displaystyle    -\theta_x \theta_y (\theta_y-1) & \displaystyle    \frac{\theta_y (\theta_y-1)}{2}  \left( \theta_x+\frac{1}{2} \right)
 \end{bmatrix}
 \]
   \begin{equation}\label{JccN}
+ \frac{\tilde{n}_y}{h}
 \begin{bmatrix}
\displaystyle \frac{\theta_x (\theta_x-1)}{2} \displaystyle \left(\theta_y + \frac{1}{2}\right) &   \displaystyle  \left(\theta_y + \frac{1}{2}\right) (1-\theta_x) (1+\theta_x) &
\boxed{\displaystyle   \frac{\theta_x (1+\theta_x)}{2} \left(\theta_y + \frac{1}{2}\right) } \\
\\
\displaystyle - \theta_y \theta_x (\theta_x-1) & \displaystyle -2 \theta_y (1-\theta_x) (1+\theta_x) &
\displaystyle -\theta_y \theta_x (1+\theta_x)  \\
 \\
 \displaystyle  \frac{\theta_x (\theta_x-1)}{2} \left(\theta_y - \frac{1}{2}\right) & \displaystyle   (1-\theta_x) (1+\theta_x) \left(\theta_y - \frac{1}{2}\right) &
  \displaystyle  \frac{\theta_x (1+\theta_x)}{2}  \left(\theta_y - \frac{1}{2}\right)
 \end{bmatrix},
     \end{equation}
with
\[
\tilde{n}_x = \frac{ \tilde{n}^*_x }{\sqrt{(\tilde{n}^*_x)^2+(\tilde{n}^*_y)^2}}, \quad \tilde{n}_y = \frac{\tilde{n}^*_y}{\sqrt{(\tilde{n}^*_x)^2+(\tilde{n}^*_y)^2}},
\]
\[
\tilde{n}^*_x =
\frac{1}{h}
 \begin{bmatrix}
\varphi_{i-2,j} & \varphi_{i-1,j} &\varphi_{i,j}  \\
 \varphi_{i-2,j-1} & \varphi_{i-1,j-1} &\varphi_{i,j-1}  \\
\varphi_{i-2,j-2} & \varphi_{i-1,j-2} &\varphi_{i,j-2}  
 \end{bmatrix}
 \odot
  \begin{bmatrix}
\displaystyle \frac{\theta_y (1+\theta_y)}{2} \displaystyle \left(\theta_x - \frac{1}{2}\right) &   \displaystyle  -\theta_x\theta_y (1+\theta_y) &  \displaystyle   \frac{\theta_y (1+\theta_y)}{2} \left( \theta_x + \frac{1}{2} \right)  \\
\\
\displaystyle (1-\theta_y) (1+\theta_y) \left(\theta_x - \frac{1}{2}\right) & \displaystyle -2 \theta_x  (1-\theta_y) (1+\theta_y) & \displaystyle (1-\theta_y) (1+\theta_y) \left( \theta_x+\frac{1}{2} \right)  \\
 \\
 \displaystyle  \frac{\theta_y (\theta_y-1)}{2} \left(\theta_x - \frac{1}{2}\right) & \displaystyle    -\theta_x \theta_y (\theta_y-1) & \displaystyle    \frac{\theta_y (\theta_y-1)}{2}  \left( \theta_x+\frac{1}{2} \right)
 \end{bmatrix},
\]
\[
\tilde{n}^*_y =
\frac{1}{h}
 \begin{bmatrix}
\varphi_{i-2,j} & \varphi_{i-1,j} &\varphi_{i,j}  \\
 \varphi_{i-2,j-1} & \varphi_{i-1,j-1} &\varphi_{i,j-1}  \\
\varphi_{i-2,j-2} & \varphi_{i-1,j-2} &\varphi_{i,j-2}  
 \end{bmatrix}
 \odot
 \begin{bmatrix}
\displaystyle \frac{\theta_x (\theta_x-1)}{2} \displaystyle \left(\theta_y + \frac{1}{2}\right) &   \displaystyle  \left(\theta_y + \frac{1}{2}\right) (1-\theta_x) (1+\theta_x) &
\displaystyle   \frac{\theta_x (1+\theta_x)}{2} \left(\theta_y + \frac{1}{2}\right) \\
\\
\displaystyle - \theta_y \theta_x (\theta_x-1) & \displaystyle -2 \theta_y (1-\theta_x) (1+\theta_x) &
\displaystyle -\theta_y \theta_x (1+\theta_x)  \\
 \\
 \displaystyle  \frac{\theta_x (\theta_x-1)}{2} \left(\theta_y - \frac{1}{2}\right) & \displaystyle   (1-\theta_x) (1+\theta_x) \left(\theta_y - \frac{1}{2}\right) &
  \displaystyle  \frac{\theta_x (1+\theta_x)}{2}  \left(\theta_y - \frac{1}{2}\right)
 \end{bmatrix},
\]
where the product operator $\odot$ means the product component-wise between the two $3 \times 3$ matrices and the sum over all components (inner products between the two vector representations of the two $3 \time 3$ matrices).

The parameter $\tau^s$ of Eq.~\eqref{CFLcond} must satisfy the condition:
\begin{equation}\label{CFLcondN}
\left| 1 - \tau^s   \left( 
\displaystyle \frac{|\tilde{n}_x| }{h} \displaystyle   \frac{\theta_y (1+\theta_y)}{2} \left( \theta_x + \frac{1}{2} \right) +
\displaystyle \frac{|\tilde{n}_y| }{h}  \displaystyle   \frac{\theta_y (1+\theta_y)}{2} \left( \theta_x + \frac{1}{2} \right)
\right)  \right| < 1.
\end{equation}
We choose a parameter $\tau^s$ that satisfies the condition \eqref{CFLcondN} for any $0\leq \theta_x, \theta_y \leq 1$ and for any $(\tilde{n}_{x},\tilde{n}_{y}) \colon (\tilde{n}_{x})^2+(\tilde{n}_{y})^2 = 1$. This is achieved by:
$
0 < \tau^s < \displaystyle \frac{2 \sqrt{2} \, h}{3}.
$
For practical purposes, we use
$
\tau^s = 0.9 \displaystyle \frac{2 \sqrt{2} \, h}{3}.
$

\section{Numerical tests}\label{sec:numtest}
In this section we confirm numerically the second order accuracy of the numerical method and we evaluate the efficiency of the multigrid approach.
We choose the following parameters in Eq.~\eqref{maineq} (see~\cite{Matteo:monum}):
\begin{equation}\label{dataeq}
a=10^4, \quad d=0.1, \quad, m_s = 64.06, \quad m_c = 100.09, \quad, \phi(c)=0.1 + 0.01 c.
\end{equation}
We choose $\Delta t = \Delta x = \Delta y = h$ and we compute the numerical solution up to the final time $t=1$.
The W-cycle iteration scheme of the multigrid is performed with $\nu_1=2$ pre-relaxation, $\nu_2=1$ post-relaxation and with an $8 \times 8$ grid as the coarsest grid.

\subsection{Accuracy tests}\label{ssec:numtest:accuracy}
In order to test the accuracy, we modify the numerical method to solve a more general problem than~\eqref{maineq} by adding source terms $f_1, f_2 \colon \Omega \times [0,T] \rightarrow \R$: 
\begin{equation*}
\left\{
\begin{array}{rcll}
\displaystyle \frac{\partial \left( \phi(c)s \right)}{\partial t} &=& \displaystyle -\frac{a}{m_c} \phi(c) \, s \, c+ d\, \nabla \cdot (\phi(c) \nabla s) + f_1 & \mbox{ in } \Omega \times [0,T]\\
\displaystyle \frac{\partial c}{\partial t} &=& \displaystyle -\frac{a}{m_s} \phi (c) \, s\, c + f_2 & \mbox{ in } \Omega \times [0,T] \\
s(\vec{x},t) &=& s_b(\vec{x},t) & \mbox{ for } (\vec{x},t) \in \partial \Omega \times [0,T] \\
s(\vec{x},0) &=& s_0(\vec{x}) & \mbox{ for } \vec{x} \in \Omega \\
c(\vec{x},0) &=& c_0(\vec{x}) & \mbox{ for } \vec{x} \in \Omega
\end{array}
\right.
\end{equation*}

We choose $f_1,f_2,s_b$ in such a way that the exact solutions are:
\[
s^\text{exa} (\vec{x},t) = 2+\sin(x) \cos(y) \sin(t+\sqrt{2}), \quad 
c^\text{exa} (\vec{x},t) = 3+\sin(0.5 \, x) \cos(3 \, y) \sin(2 \, t+\sqrt{3})
\]
and then we compute the $L^p$ errors at time $t=1$ on the solutions
\[
e^s_h = \frac{\left\| s_h -s_h^\text{exa} \right\|_p}{\left\| s_h^\text{exa} \right\|_p}, \quad 
e^c_h = \frac{\left\| c_h -c_h^\text{exa} \right\|_p}{\left\| c_h^\text{exa} \right\|_p}
\]
and on the gradients 
\[
e^{|\nabla s|}_h = \frac{\left\| |\nabla s_h| -|\nabla s_h^\text{exa}| \right\|_p}{\left\| | \nabla s_h^\text{exa}| \right\|_p}, \quad 
e^{|\nabla c|}_h = \frac{\left\| |\nabla c_h| -|\nabla c_h^\text{exa}| \right\|_p}{\left\| | \nabla c_h^\text{exa}| \right\|_p},
\]
where $s_h$ and $c_h$ are the numerical solutions and $\nabla s_h$ and $\nabla c_h$ are computed by central differences.
We perform two tests. In {\sc Test 1}, the domain is represented by a circle and the level-set function is:
\begin{equation}\label{test13}
\varphi(x,y) = \sqrt{(x-x_0)^2+(y-y_0)^2}-R, \quad \text{ where } x_0 = \sqrt{2}/30, \; y_0 = \sqrt{3}/40, \; R = 1.486.
\end{equation}

\begin{figure}[!hbt]
  \begin{minipage}[c]{0.45\textwidth}
   	\centering
   	\includegraphics[width=0.79\textwidth]{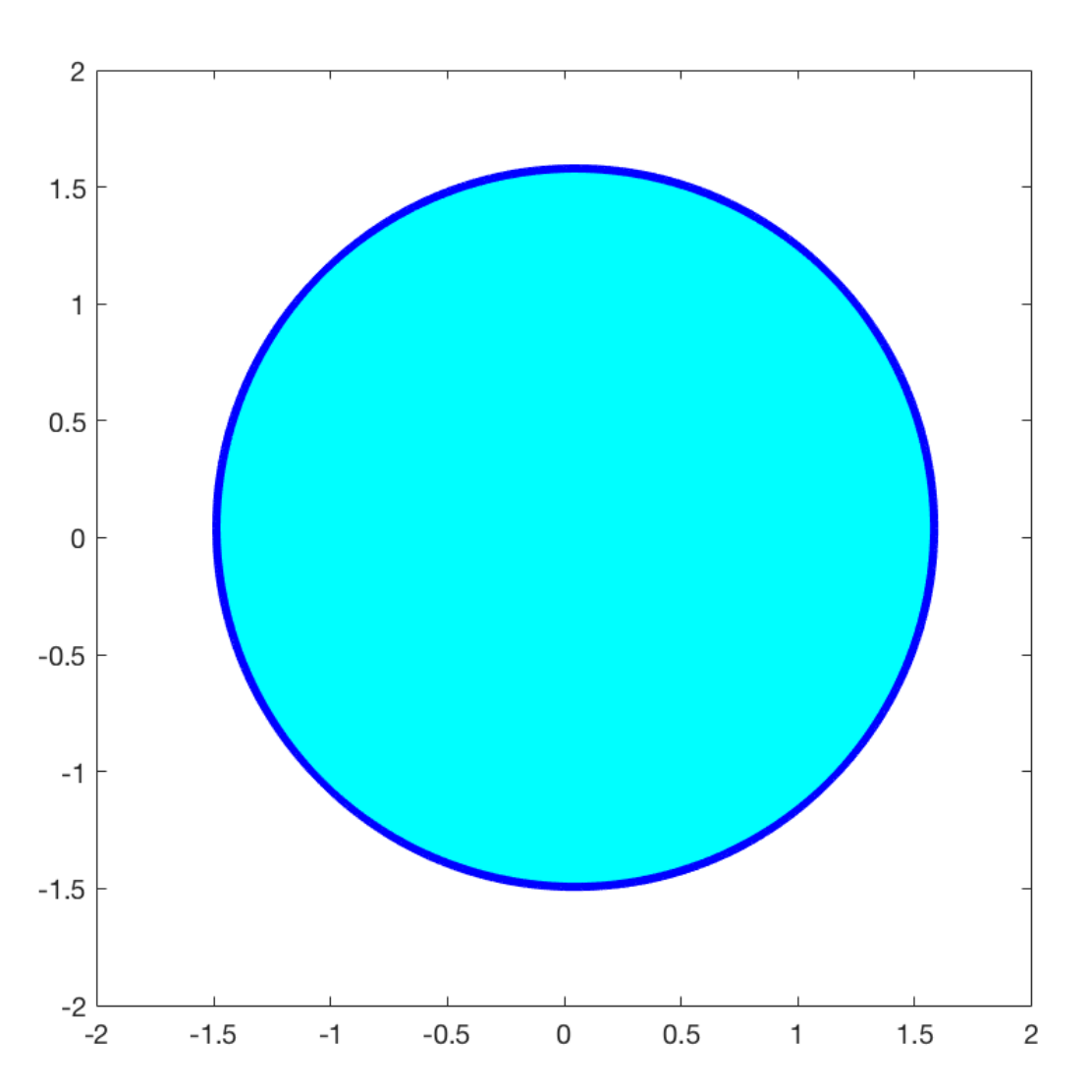}
 \captionsetup{width=0.80\textwidth}
 \end{minipage}
 \ \hspace{2mm} \hspace{3mm} \
 \begin{minipage}[c]{0.45\textwidth}
   	\centering
   	\includegraphics[width=0.79\textwidth]{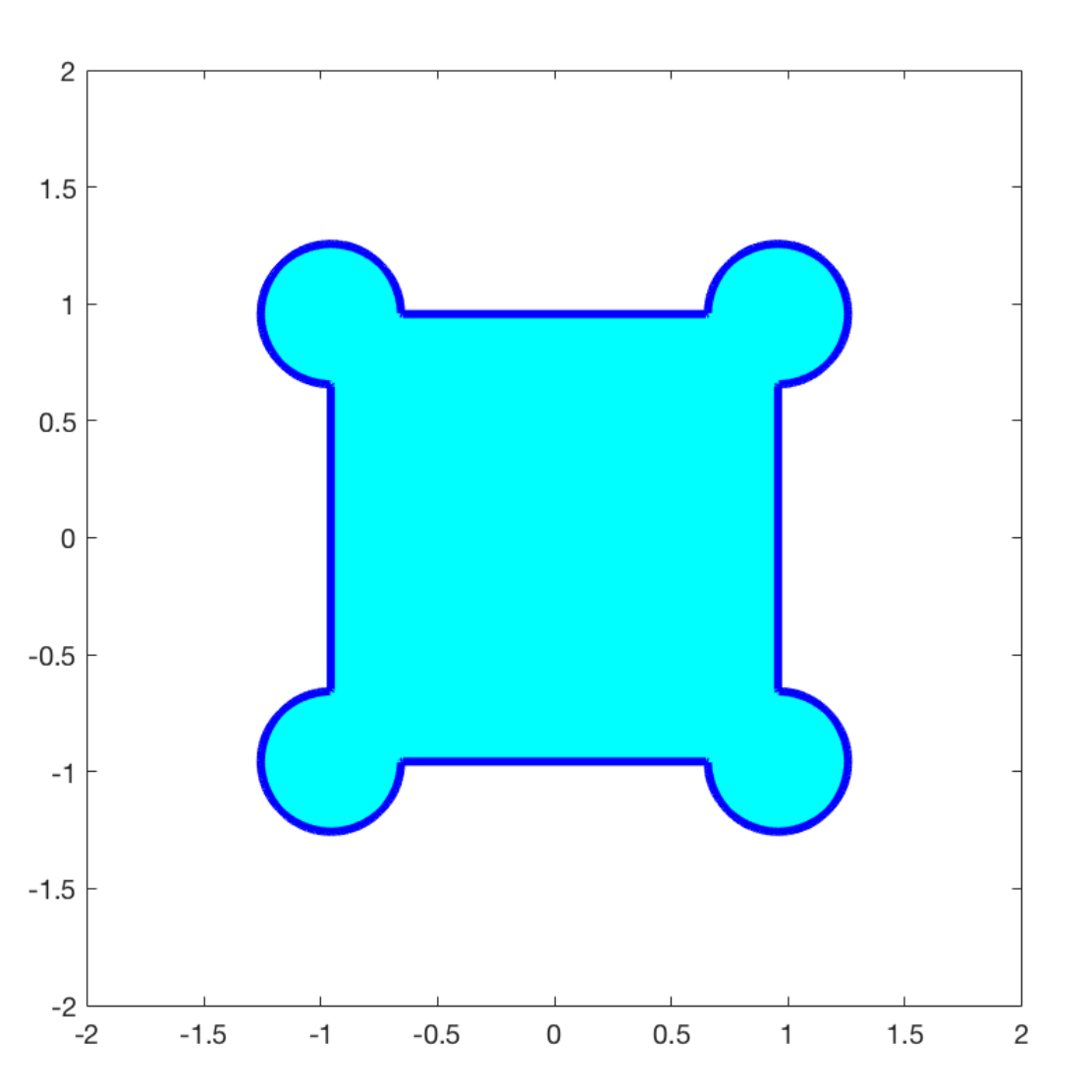}
 \captionsetup{width=0.80\textwidth}
	 \end{minipage}
	\caption{{ Domain $\Omega$ for {\sc Test 1, Test 1N and 3} (left) and {\sc Test 2, Test 2N and 4} (right).}}
	\label{fig:domains}
\end{figure}

\begin{figure}[!hbt]
 \begin{minipage}[c]{0.45\textwidth}
   	\centering
   	\includegraphics[width=0.79\textwidth]{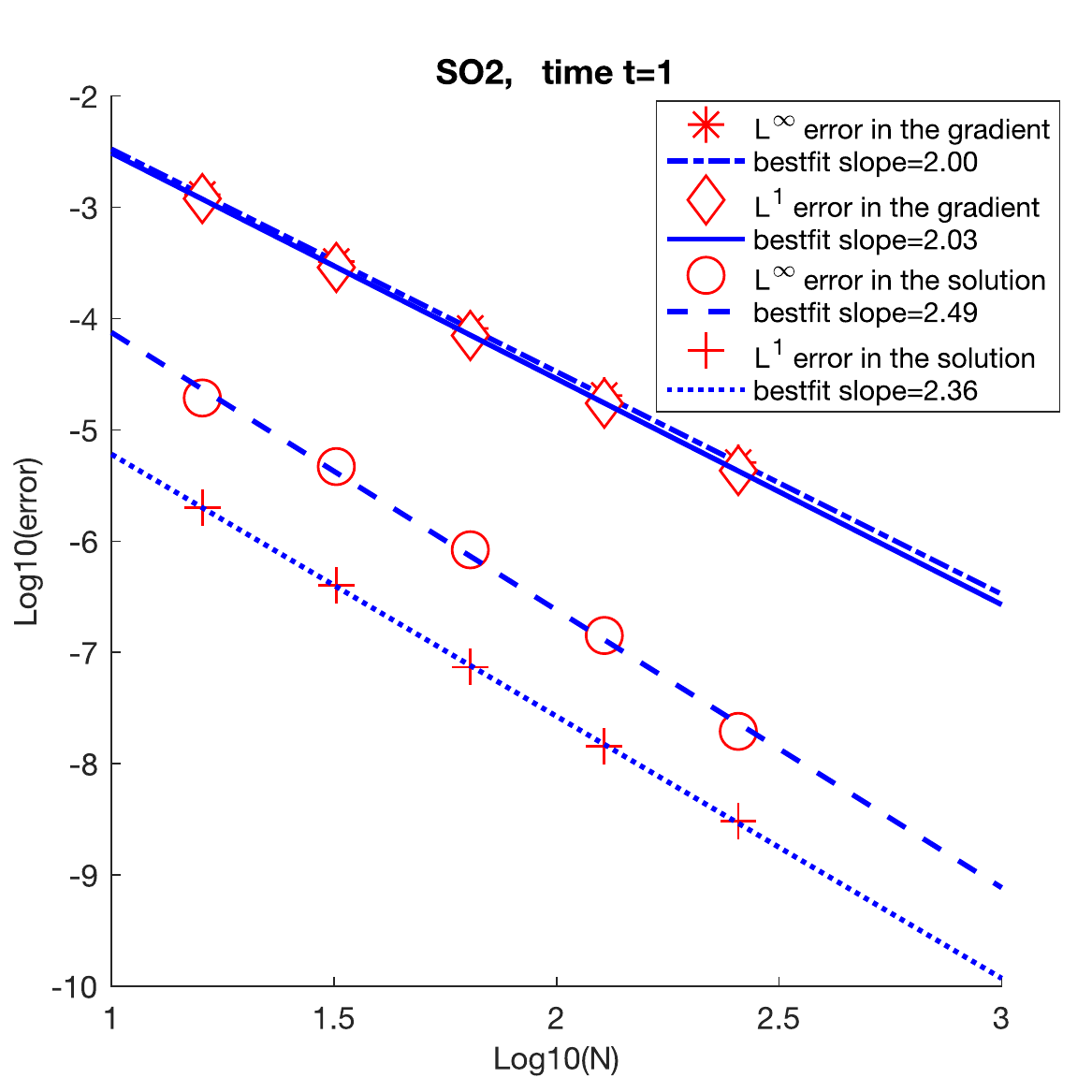}
 \captionsetup{width=0.80\textwidth}
	\caption{{ {\sc Test 1}: bestfit lines of the errors in the solution and in the gradient for $\text{SO}_2$ in $L^1$ and $L^\infty$ norms (Table \ref{table:test1s}).}}
	\label{fig:test1s}
 \end{minipage}
 \ \hspace{2mm} \hspace{3mm} \
 \begin{minipage}[c]{0.45\textwidth}
\centering
   	\includegraphics[width=0.79\textwidth]{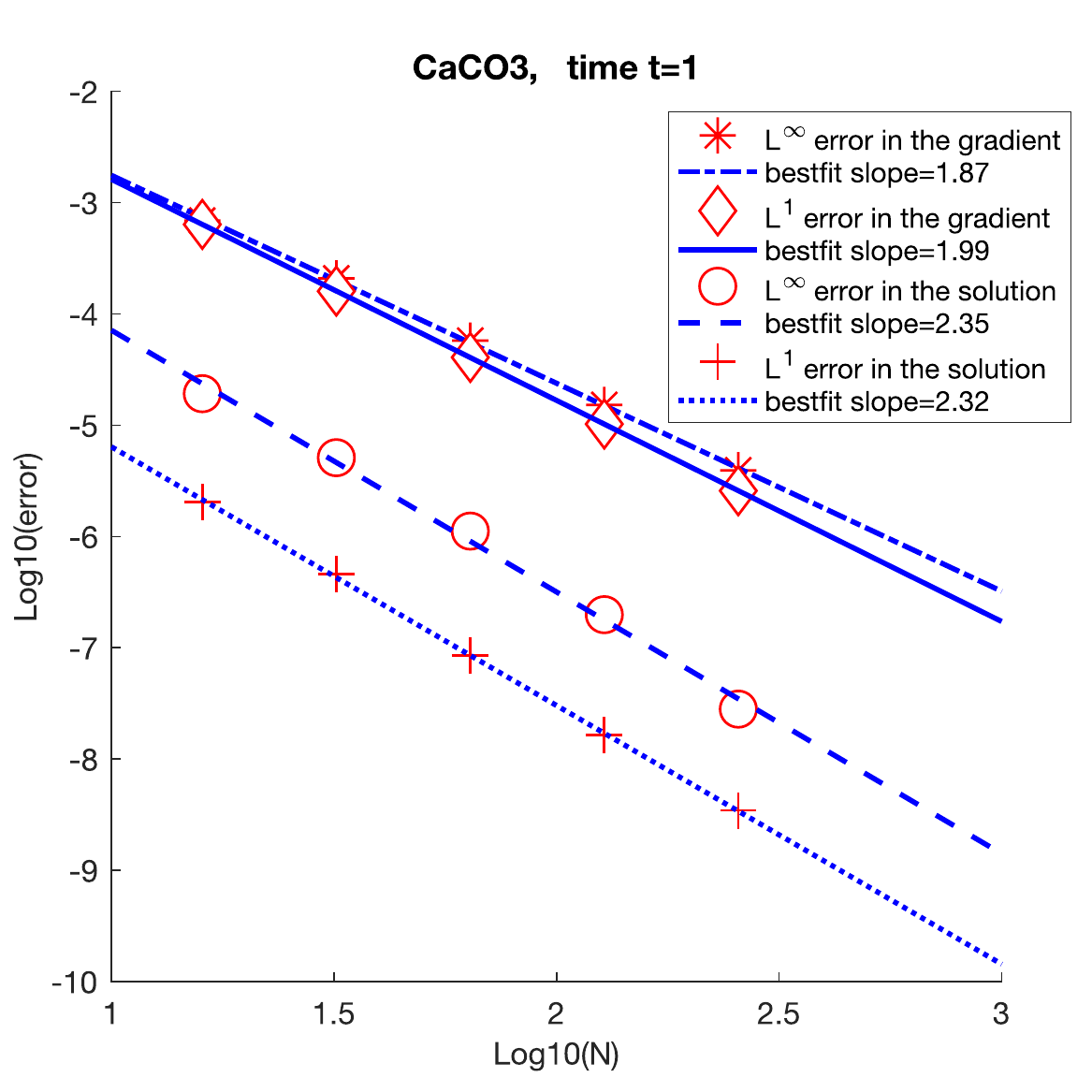}
 \captionsetup{width=0.80\textwidth}
	\caption{{ {\sc Test 1}: bestfit lines of the errors in the solution and in the gradient for $\text{CaCO}_3$ in $L^1$ and $L^\infty$ norms (Table \ref{table:test1c}).}}
	\label{fig:test1c}
 \end{minipage}
\end{figure}

\begin{figure}[!hbt]
  \begin{minipage}[c]{0.45\textwidth}
   	\centering
   	\includegraphics[width=0.79\textwidth]{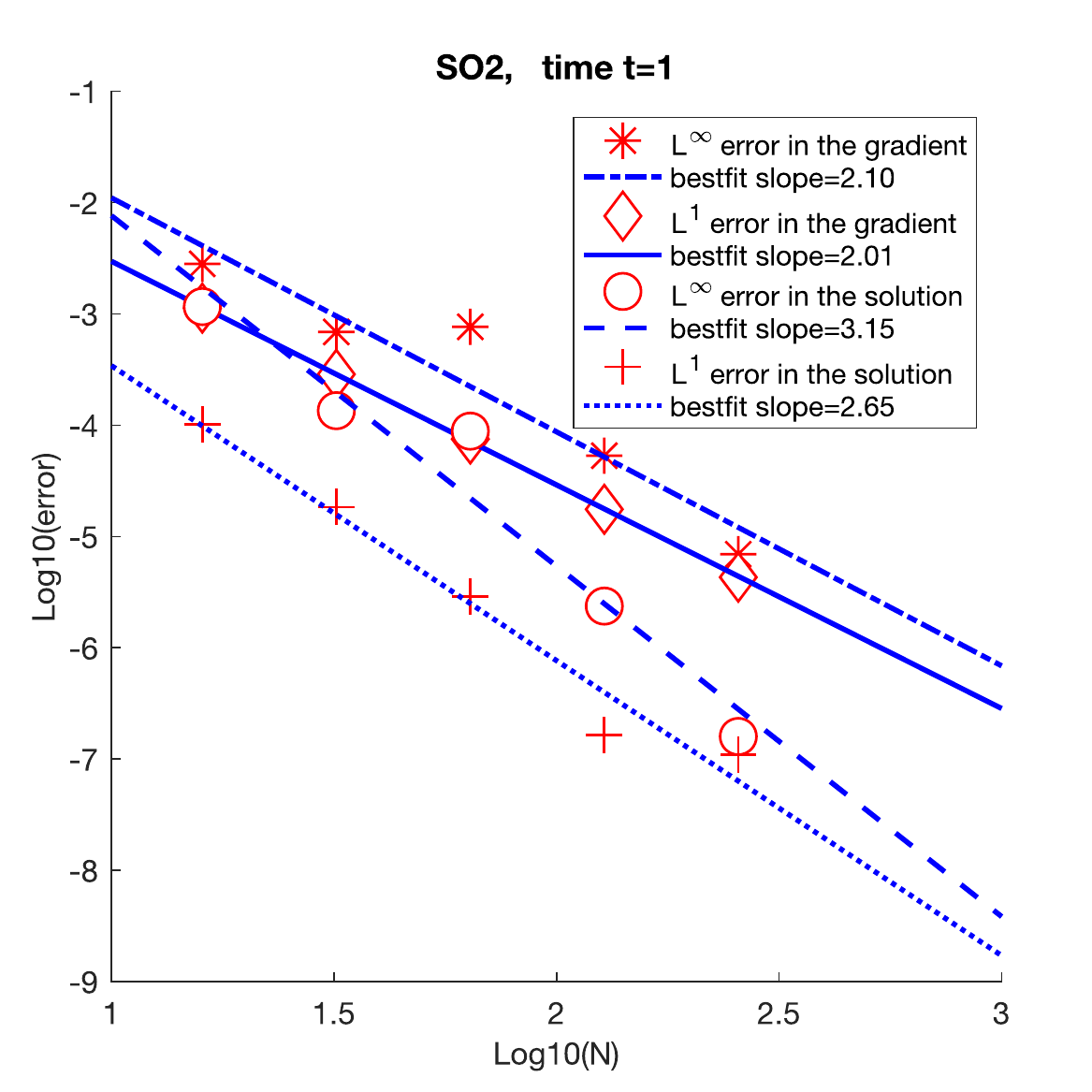}
 \captionsetup{width=0.80\textwidth}
	\caption{{ {\sc Test 2}: bestfit lines of the errors in the solution and in the gradient for $\text{SO}_2$ in $L^1$ and $L^\infty$ norms (Table \ref{table:test2s}).}}
	\label{fig:test2s}
 \end{minipage}
 \ \hspace{2mm} \hspace{3mm} \
 \begin{minipage}[c]{0.45\textwidth}
   	\centering
   	\includegraphics[width=0.79\textwidth]{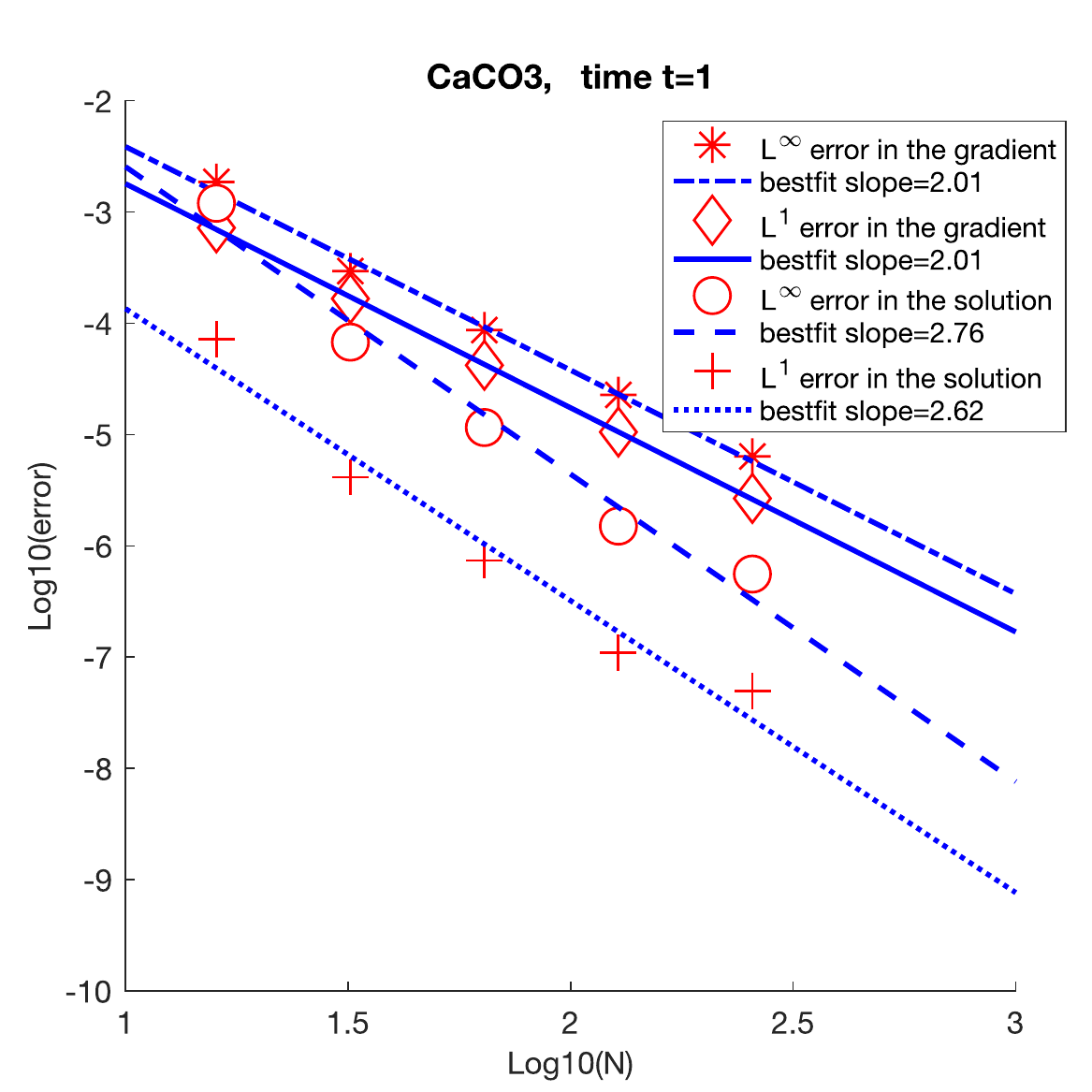}
 \captionsetup{width=0.80\textwidth}
	\caption{{ {\sc Test 2}: bestfit lines of the errors in the solution and in the gradient for $\text{CaCO}_3$ in $L^1$ and $L^\infty$ norms (Table \ref{table:test2c}).}}
	\label{fig:test2c}
\end{minipage}
\end{figure}

\begin{figure}[!hbt]
 \begin{minipage}[c]{0.45\textwidth}
   	\centering
   	\includegraphics[width=0.79\textwidth]{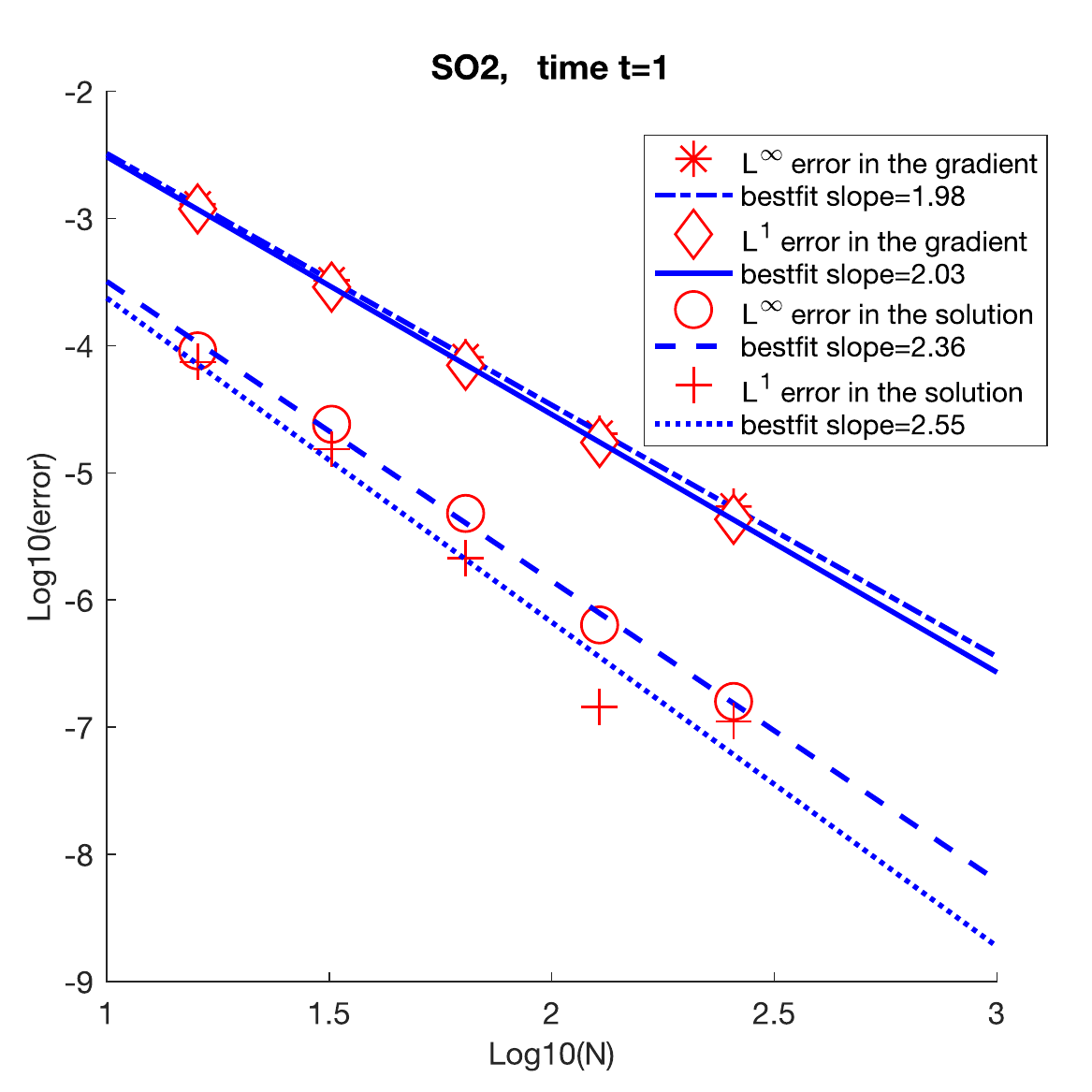}
 \captionsetup{width=0.80\textwidth}
	\caption{{ {\sc Test 1N}: bestfit lines of the errors in the solution and in the gradient for $\text{SO}_2$ in $L^1$ and $L^\infty$ norms (Table \ref{table:test1sN}).}}
	\label{fig:test1sN}
 \end{minipage}
 \ \hspace{2mm} \hspace{3mm} \
 \begin{minipage}[c]{0.45\textwidth}
\centering
   	\includegraphics[width=0.79\textwidth]{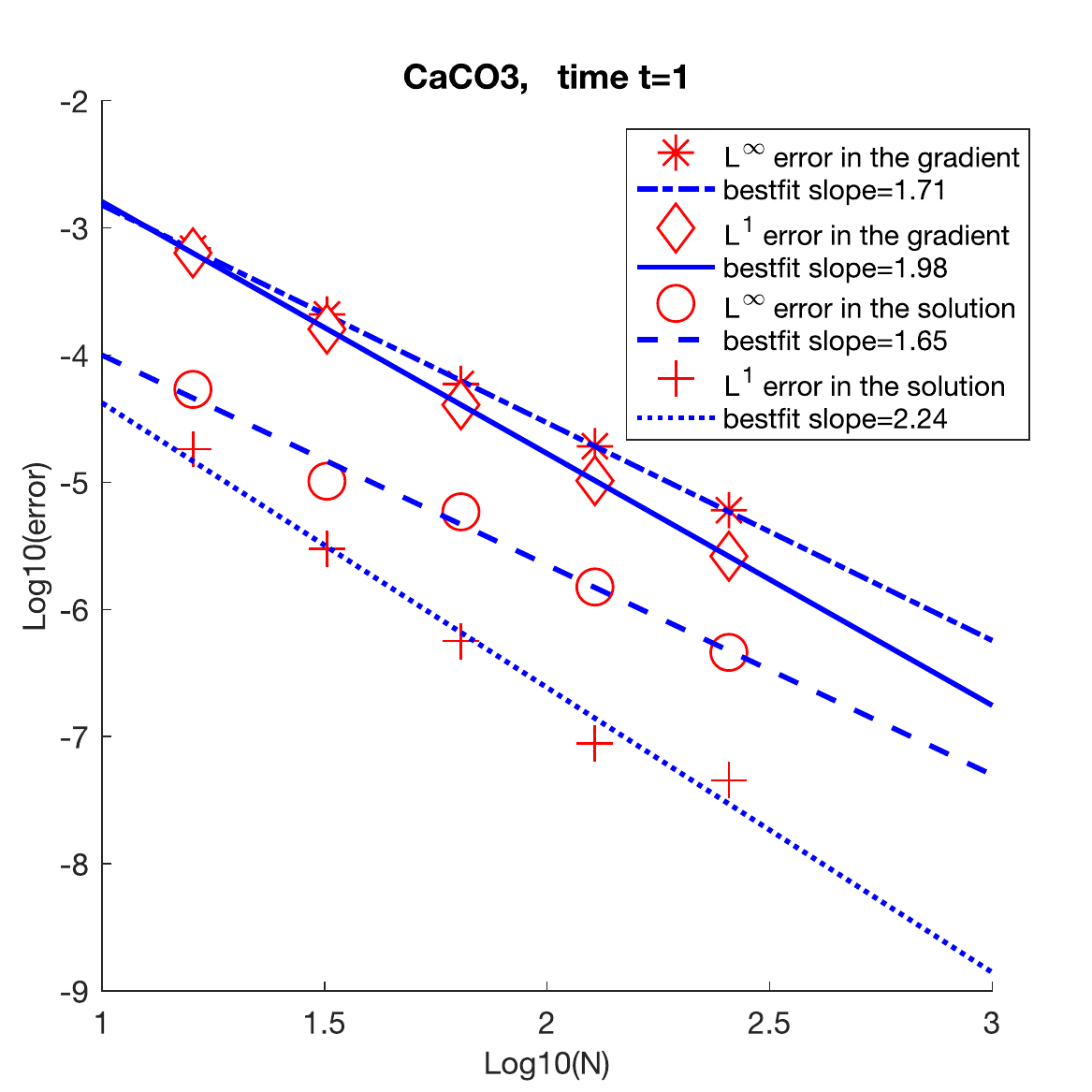}
 \captionsetup{width=0.80\textwidth}
	\caption{{ {\sc Test 1N}: bestfit lines of the errors in the solution and in the gradient for $\text{CaCO}_3$ in $L^1$ and $L^\infty$ norms (Table \ref{table:test1cN}).}}
	\label{fig:test1cN}
 \end{minipage}
\end{figure}

\begin{figure}[!hbt]
  \begin{minipage}[c]{0.45\textwidth}
   	\centering
   	\includegraphics[width=0.79\textwidth]{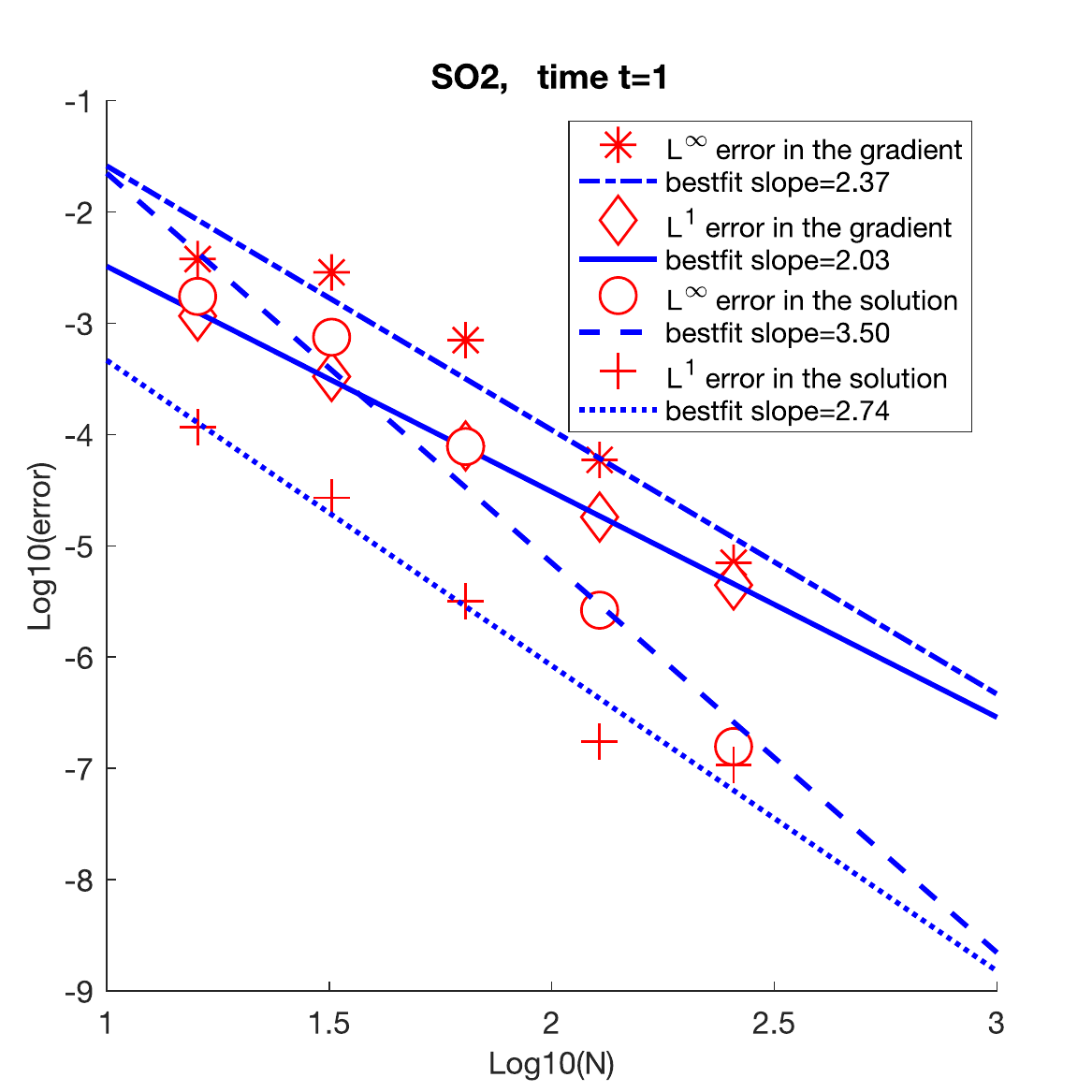}
 \captionsetup{width=0.80\textwidth}
	\caption{{ {\sc Test 2N}: bestfit lines of the errors in the solution and in the gradient for $\text{SO}_2$ in $L^1$ and $L^\infty$ norms (Table \ref{table:test2sN}).}}
	\label{fig:test2sN}
 \end{minipage}
 \ \hspace{2mm} \hspace{3mm} \
 \begin{minipage}[c]{0.45\textwidth}
   	\centering
   	\includegraphics[width=0.79\textwidth]{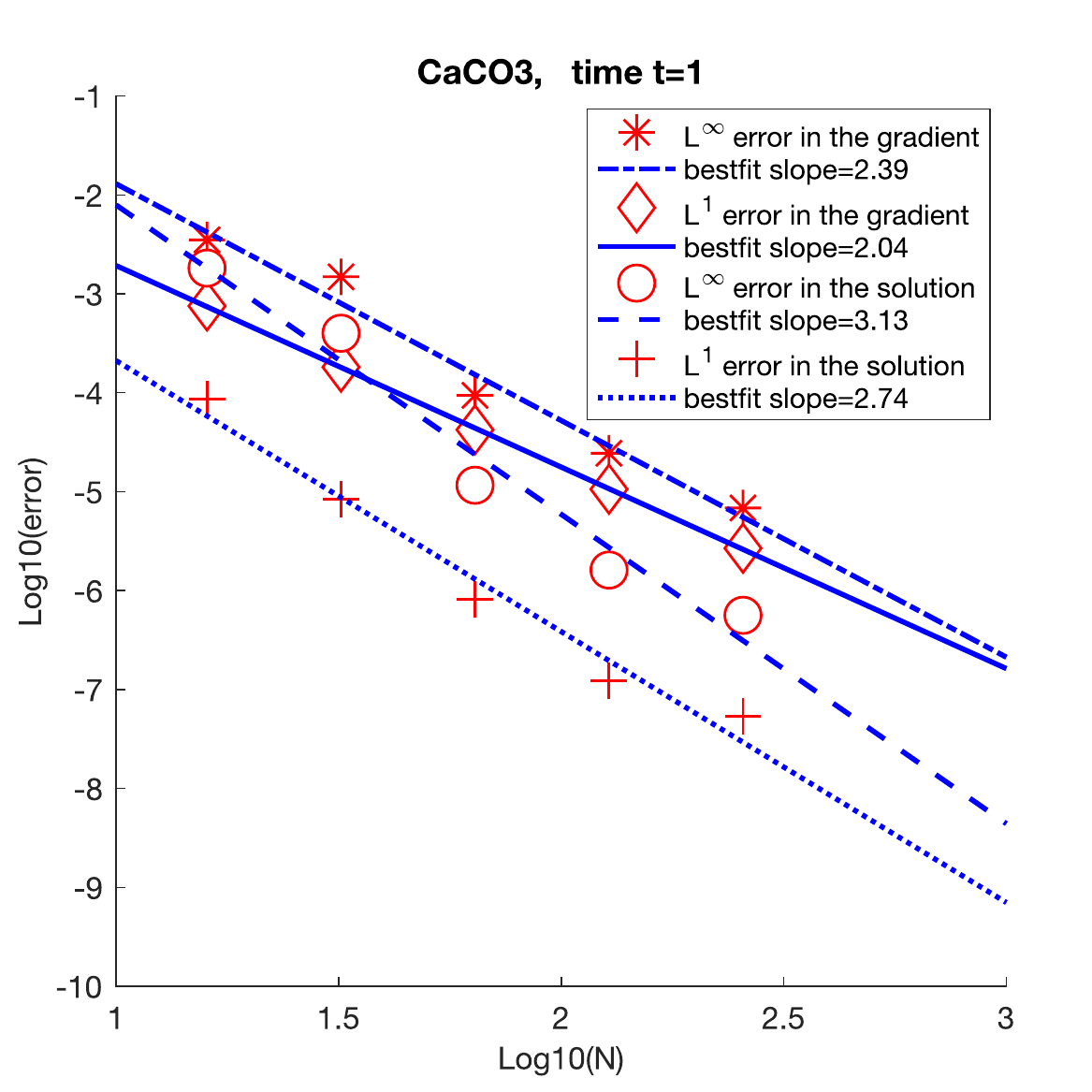}
 \captionsetup{width=0.80\textwidth}
	\caption{{ {\sc Test 2N}: bestfit lines of the errors in the solution and in the gradient for $\text{CaCO}_3$ in $L^1$ and $L^\infty$ norms (Table \ref{table:test2cN}).}}
	\label{fig:test2cN}
\end{minipage}
\end{figure}

In {\sc Test 2} the domain is represented by the union of a square and four circles centred on the vertices

\begin{equation}\label{test24}
\varphi(x,y)=\min \left\{ \varphi_1(x,y), \varphi_2(x,y) \right\}
\end{equation}
where 
\[
\varphi_1(x,y)=\max \left\{ |x|, |y| \right\} - L, \quad
\varphi_2(x,y)=\sqrt{(|x|-L)^2+(|y|-L)^2}-D, \quad
L= 0.9567, \quad D = 0.3.
\]
The domains for {\sc Test 1} and {\sc Test 2} are represented in Fig.~\ref{fig:domains}.
$L^1$ and $L^\infty$ errors for the solutions and the gradients for {\sc Test 1} can be found in Table \ref{table:test1s} (for $\text{SO}_2$) and Table \ref{table:test1c} (for $\text{CaCO}_3$), and for {\sc Test 2} in Table \ref{table:test2s} (for $\text{SO}_2$) and Table \ref{table:test2c} (for $\text{CaCO}_3$).
Bestfit lines in bilogarithmic plots for the errors versus $N$ can be found in Figs.~\ref{fig:test1s} ($\text{SO}_2$ in {\sc Test 1}), \ref{fig:test1c} ($\text{CaCO}_3$ in {\sc Test 1}), \ref{fig:test2s} ($\text{SO}_2$ in {\sc Test 2}) and \ref{fig:test2c} ($\text{CaCO}_3$ in {\sc Test 2}).
We repeat the two tests for the case of Neumann boundary conditions for $s$ as described in \S\ref{ssec:numtest:neumann}) ({\sc Test 1N} and {\sc Test 2N}) and present the results
in Tables \ref{table:test1sN}, \ref{table:test1cN}, \ref{table:test2sN} and \ref{table:test2cN} and in Figures~\ref{fig:test1sN}, \ref{fig:test1cN}, \ref{fig:test2sN} and \ref{fig:test2cN}. 
We note that in all cases second order convergence is achieved.

\begin{table}[!hbt]
\captionsetup{width=0.80\textwidth}
\caption{ { {\sc Test 1}. Accuracy order in the solution (top) and in the gradient (bottom) for $s$ ($\text{SO}_2$). 
}} 
\centering      
\begin{tabular}{|| c || c | c || c | c ||} 
\hline \hline 
No.~of points & $L^1$ error of $s$ & order & $L^\infty$ error of $s$ & order \\ 
\hline 
16 $\times$ 16 & 2.00 $\cdot 10^{-6}$ & - & 1.93 $\cdot 10^{-5}$ & - \\ 
32 $\times$ 32 & 4.03 $\cdot 10^{-7}$ & 2.31 & 4.67 $\cdot 10^{-6}$ & 2.05 \\ 
64 $\times$ 64 & 7.38 $\cdot 10^{-8}$ & 2.45 & 8.35 $\cdot 10^{-7}$ & 2.48 \\ 
128 $\times$ 128 & 1.43 $\cdot 10^{-8}$ & 2.37 & 1.42 $\cdot 10^{-7}$ & 2.56 \\ 
256 $\times$ 256 & 3.03 $\cdot 10^{-9}$ & 2.24 & 1.95 $\cdot 10^{-8}$ & 2.86 \\ 
\hline \hline 
No.~of points & $L^1$ error of $|\nabla s|$ & order & $L^\infty$ error of $|\nabla s|$ & order \\ 
\hline 
16 $\times$ 16 & 1.20 $\cdot 10^{-3}$ & - & 1.30 $\cdot 10^{-3}$ & - \\ 
32 $\times$ 32 & 2.88 $\cdot 10^{-4}$ & 2.06 & 3.25 $\cdot 10^{-4}$ & 2.00 \\ 
64 $\times$ 64 & 7.03 $\cdot 10^{-5}$ & 2.03 & 8.13 $\cdot 10^{-5}$ & 2.00 \\ 
128 $\times$ 128 & 1.73 $\cdot 10^{-5}$ & 2.02 & 2.03 $\cdot 10^{-5}$ & 2.00 \\ 
256 $\times$ 256 & 4.30 $\cdot 10^{-6}$ & 2.01 & 5.08 $\cdot 10^{-6}$ & 2.00 \\ 
\hline \hline 
\end{tabular} 
    \label{table:test1s}  
 \end{table}
 
\begin{table}[!hbt]
\captionsetup{width=0.80\textwidth}
\caption{ { {\sc Test 1}. Accuracy order in the solution (top) and in the gradient (bottom) for $c$ ($\text{CaCO}_3$). 
}} 
\centering      
\begin{tabular}{|| c || c | c || c | c ||} 
\hline \hline 
No.~of points & $L^1$ error of $c$ & order & $L^\infty$ error of $c$ & order \\ 
\hline 
16 $\times$ 16 & 2.04 $\cdot 10^{-6}$ & - & 1.92 $\cdot 10^{-5}$ & - \\ 
32 $\times$ 32 & 4.60 $\cdot 10^{-7}$ & 2.15 & 5.10 $\cdot 10^{-6}$ & 1.91 \\ 
64 $\times$ 64 & 8.52 $\cdot 10^{-8}$ & 2.43 & 1.12 $\cdot 10^{-6}$ & 2.19 \\ 
128 $\times$ 128 & 1.63 $\cdot 10^{-8}$ & 2.39 & 1.98 $\cdot 10^{-7}$ & 2.50 \\ 
256 $\times$ 256 & 3.45 $\cdot 10^{-9}$ & 2.24 & 2.80 $\cdot 10^{-8}$ & 2.82 \\ 
\hline \hline 
No.~of points & $L^1$ error of $|\nabla c|$ & order & $L^\infty$ error of $|\nabla c|$ & order \\ 
\hline 
16 $\times$ 16 & 6.36 $\cdot 10^{-4}$ & - & 6.89 $\cdot 10^{-4}$ & - \\ 
32 $\times$ 32 & 1.60 $\cdot 10^{-4}$ & 1.99 & 2.08 $\cdot 10^{-4}$ & 1.73 \\ 
64 $\times$ 64 & 4.06 $\cdot 10^{-5}$ & 1.98 & 5.74 $\cdot 10^{-5}$ & 1.86 \\ 
128 $\times$ 128 & 1.02 $\cdot 10^{-5}$ & 1.99 & 1.52 $\cdot 10^{-5}$ & 1.92 \\ 
256 $\times$ 256 & 2.57 $\cdot 10^{-6}$ & 1.99 & 3.92 $\cdot 10^{-6}$ & 1.95 \\ 
\hline \hline 
\end{tabular} 
    \label{table:test1c}  
 \end{table}
 
 \begin{table}[!hbt]
\captionsetup{width=0.80\textwidth}
\caption{ { {\sc Test 2}. Accuracy order in the solution (top) and in the gradient (bottom) for $s$ ($\text{SO}_2$). 
}} 
\centering      
\begin{tabular}{|| c || c | c || c | c ||} 
\hline \hline 
No.~of points & $L^1$ error of $s$ & order & $L^\infty$ error of $s$ & order \\ 
\hline 
16 $\times$ 16 & 1.01 $\cdot 10^{-4}$ & - & 1.16 $\cdot 10^{-3}$ & - \\ 
32 $\times$ 32 & 1.83 $\cdot 10^{-5}$ & 2.47 & 1.35 $\cdot 10^{-4}$ & 3.11 \\ 
64 $\times$ 64 & 2.88 $\cdot 10^{-6}$ & 2.67 & 8.83 $\cdot 10^{-5}$ & 0.61 \\ 
128 $\times$ 128 & 1.65 $\cdot 10^{-7}$ & 4.13 & 2.36 $\cdot 10^{-6}$ & 5.22 \\ 
256 $\times$ 256 & 1.10 $\cdot 10^{-7}$ & 0.59 & 1.59 $\cdot 10^{-7}$ & 3.89 \\ 
\hline \hline 
No.~of points & $L^1$ error of $|\nabla s|$ & order & $L^\infty$ error of $|\nabla s|$ & order \\ 
\hline 
16 $\times$ 16 & 1.12 $\cdot 10^{-3}$ & - & 2.81 $\cdot 10^{-3}$ & - \\ 
32 $\times$ 32 & 2.86 $\cdot 10^{-4}$ & 1.97 & 6.87 $\cdot 10^{-4}$ & 2.03 \\ 
64 $\times$ 64 & 7.49 $\cdot 10^{-5}$ & 1.93 & 7.64 $\cdot 10^{-4}$ & -0.15 \\ 
128 $\times$ 128 & 1.75 $\cdot 10^{-5}$ & 2.10 & 5.31 $\cdot 10^{-5}$ & 3.85 \\ 
256 $\times$ 256 & 4.30 $\cdot 10^{-6}$ & 2.03 & 6.92 $\cdot 10^{-6}$ & 2.94 \\ 
\hline \hline 
\end{tabular} 
    \label{table:test2s}  
 \end{table}
 
\begin{table}[!hbt]
\captionsetup{width=0.80\textwidth}
\caption{ { {\sc Test 2}. Accuracy order in the solution (top) and in the gradient (bottom) for $c$ ($\text{CaCO}_3$). 
}} 
\centering      
\begin{tabular}{|| c || c | c || c | c ||} 
\hline \hline 
No.~of points & $L^1$ error of $c$ & order & $L^\infty$ error of $c$ & order \\ 
\hline 
16 $\times$ 16 & 7.20 $\cdot 10^{-5}$ & - & 1.20 $\cdot 10^{-3}$ & - \\ 
32 $\times$ 32 & 4.13 $\cdot 10^{-6}$ & 4.12 & 6.76 $\cdot 10^{-5}$ & 4.15 \\ 
64 $\times$ 64 & 7.37 $\cdot 10^{-7}$ & 2.49 & 1.15 $\cdot 10^{-5}$ & 2.55 \\ 
128 $\times$ 128 & 1.10 $\cdot 10^{-7}$ & 2.74 & 1.50 $\cdot 10^{-6}$ & 2.94 \\ 
256 $\times$ 256 & 4.96 $\cdot 10^{-8}$ & 1.15 & 5.57 $\cdot 10^{-7}$ & 1.43 \\ 
\hline \hline 
No.~of points & $L^1$ error of $|\nabla c|$ & order & $L^\infty$ error of $|\nabla c|$ & order \\ 
\hline 
16 $\times$ 16 & 7.20 $\cdot 10^{-4}$ & - & 1.86 $\cdot 10^{-3}$ & - \\ 
32 $\times$ 32 & 1.66 $\cdot 10^{-4}$ & 2.12 & 2.94 $\cdot 10^{-4}$ & 2.66 \\ 
64 $\times$ 64 & 4.19 $\cdot 10^{-5}$ & 1.99 & 8.68 $\cdot 10^{-5}$ & 1.76 \\ 
128 $\times$ 128 & 1.05 $\cdot 10^{-5}$ & 1.99 & 2.27 $\cdot 10^{-5}$ & 1.94 \\ 
256 $\times$ 256 & 2.67 $\cdot 10^{-6}$ & 1.98 & 6.35 $\cdot 10^{-6}$ & 1.84 \\ 
\hline \hline 
\end{tabular} 
    \label{table:test2c}  
 \end{table}

\begin{table}[!hbt]
\captionsetup{width=0.80\textwidth}
\caption{ { {\sc Test 1N}. Accuracy order in the solution (top) and in the gradient (bottom) for $s$ ($\text{SO}_2$). 
}} 
\centering      
\begin{tabular}{|| c || c | c || c | c ||} 
\hline \hline 
No.~of points & $L^1$ error of $s$ & order & $L^\infty$ error of $s$ & order \\ 
\hline 
16 $\times$ 16 & 7.43 $\cdot 10^{-5}$ & - & 9.10 $\cdot 10^{-5}$ & - \\ 
32 $\times$ 32 & 1.53 $\cdot 10^{-5}$ & 2.28 & 2.41 $\cdot 10^{-5}$ & 1.92 \\ 
64 $\times$ 64 & 2.14 $\cdot 10^{-6}$ & 2.84 & 4.79 $\cdot 10^{-6}$ & 2.33 \\ 
128 $\times$ 128 & 1.45 $\cdot 10^{-7}$ & 3.88 & 6.37 $\cdot 10^{-7}$ & 2.91 \\ 
256 $\times$ 256 & 1.11 $\cdot 10^{-7}$ & 0.39 & 1.59 $\cdot 10^{-7}$ & 2.00 \\ 
\hline \hline 
No.~of points & $L^1$ error of $|\nabla s|$ & order & $L^\infty$ error of $|\nabla s|$ & order \\ 
\hline 
16 $\times$ 16 & 1.19 $\cdot 10^{-3}$ & - & 1.29 $\cdot 10^{-3}$ & - \\ 
32 $\times$ 32 & 2.89 $\cdot 10^{-4}$ & 2.04 & 3.26 $\cdot 10^{-4}$ & 1.98 \\ 
64 $\times$ 64 & 7.01 $\cdot 10^{-5}$ & 2.04 & 8.11 $\cdot 10^{-5}$ & 2.01 \\ 
128 $\times$ 128 & 1.73 $\cdot 10^{-5}$ & 2.02 & 2.03 $\cdot 10^{-5}$ & 2.00 \\ 
256 $\times$ 256 & 4.31 $\cdot 10^{-6}$ & 2.01 & 5.44 $\cdot 10^{-6}$ & 1.90 \\ 
\hline \hline 
\end{tabular} 
    \label{table:test1sN}  
 \end{table}
 
\begin{table}[!hbt]
\captionsetup{width=0.80\textwidth}
\caption{ { {\sc Test 1N}. Accuracy order in the solution (top) and in the gradient (bottom) for $c$ ($\text{CaCO}_3$). 
}} 
\centering      
\begin{tabular}{|| c || c | c || c | c ||} 
\hline \hline 
No.~of points & $L^1$ error of $c$ & order & $L^\infty$ error of $c$ & order \\ 
\hline 
16 $\times$ 16 & 1.82 $\cdot 10^{-5}$ & - & 5.37 $\cdot 10^{-5}$ & - \\ 
32 $\times$ 32 & 3.00 $\cdot 10^{-6}$ & 2.60 & 1.02 $\cdot 10^{-5}$ & 2.40 \\ 
64 $\times$ 64 & 5.63 $\cdot 10^{-7}$ & 2.41 & 5.85 $\cdot 10^{-6}$ & 0.80 \\ 
128 $\times$ 128 & 8.83 $\cdot 10^{-8}$ & 2.67 & 1.50 $\cdot 10^{-6}$ & 1.96 \\ 
256 $\times$ 256 & 4.53 $\cdot 10^{-8}$ & 0.96 & 4.59 $\cdot 10^{-7}$ & 1.71 \\ 
\hline \hline 
No.~of points & $L^1$ error of $|\nabla c|$ & order & $L^\infty$ error of $|\nabla c|$ & order \\ 
\hline 
16 $\times$ 16 & 6.35 $\cdot 10^{-4}$ & - & 6.88 $\cdot 10^{-4}$ & - \\ 
32 $\times$ 32 & 1.60 $\cdot 10^{-4}$ & 1.99 & 2.09 $\cdot 10^{-4}$ & 1.72 \\ 
64 $\times$ 64 & 4.06 $\cdot 10^{-5}$ & 1.98 & 5.92 $\cdot 10^{-5}$ & 1.82 \\ 
128 $\times$ 128 & 1.03 $\cdot 10^{-5}$ & 1.98 & 1.92 $\cdot 10^{-5}$ & 1.63 \\ 
256 $\times$ 256 & 2.62 $\cdot 10^{-6}$ & 1.97 & 6.02 $\cdot 10^{-6}$ & 1.67 \\ 
\hline \hline 
\end{tabular} 
    \label{table:test1cN}  
 \end{table}
 
 \begin{table}[!hbt]
\captionsetup{width=0.80\textwidth}
\caption{ { {\sc Test 2N}. Accuracy order in the solution (top) and in the gradient (bottom) for $s$ ($\text{SO}_2$). 
}} 
\centering      
\begin{tabular}{|| c || c | c || c | c ||} 
\hline \hline 
No.~of points & $L^1$ error of $s$ & order & $L^\infty$ error of $s$ & order \\ 
\hline 
16 $\times$ 16 & 1.16 $\cdot 10^{-4}$ & - & 1.74 $\cdot 10^{-3}$ & - \\ 
32 $\times$ 32 & 2.70 $\cdot 10^{-5}$ & 2.10 & 7.48 $\cdot 10^{-4}$ & 1.22 \\ 
64 $\times$ 64 & 3.18 $\cdot 10^{-6}$ & 3.09 & 7.79 $\cdot 10^{-5}$ & 3.26 \\ 
128 $\times$ 128 & 1.74 $\cdot 10^{-7}$ & 4.19 & 2.64 $\cdot 10^{-6}$ & 4.89 \\ 
256 $\times$ 256 & 1.07 $\cdot 10^{-7}$ & 0.69 & 1.57 $\cdot 10^{-7}$ & 4.07 \\ 
\hline \hline 
No.~of points & $L^1$ error of $|\nabla s|$ & order & $L^\infty$ error of $|\nabla s|$ & order \\ 
\hline 
16 $\times$ 16 & 1.16 $\cdot 10^{-3}$ & - & 3.77 $\cdot 10^{-3}$ & - \\ 
32 $\times$ 32 & 3.32 $\cdot 10^{-4}$ & 1.81 & 2.86 $\cdot 10^{-3}$ & 0.40 \\ 
64 $\times$ 64 & 7.89 $\cdot 10^{-5}$ & 2.07 & 7.05 $\cdot 10^{-4}$ & 2.02 \\ 
128 $\times$ 128 & 1.81 $\cdot 10^{-5}$ & 2.12 & 5.91 $\cdot 10^{-5}$ & 3.58 \\ 
256 $\times$ 256 & 4.44 $\cdot 10^{-6}$ & 2.03 & 7.03 $\cdot 10^{-6}$ & 3.07 \\ 
\hline \hline 
\end{tabular} 	
    \label{table:test2sN}  
 \end{table}
 
\begin{table}[!hbt]
\captionsetup{width=0.80\textwidth}
\caption{ { {\sc Test 2N}. Accuracy order in the solution (top) and in the gradient (bottom) for $c$ ($\text{CaCO}_3$). 
}} 
\centering      
\begin{tabular}{|| c || c | c || c | c ||} 
\hline \hline 
No.~of points & $L^1$ error of $c$ & order & $L^\infty$ error of $c$ & order \\ 
\hline 
16 $\times$ 16 & 8.57 $\cdot 10^{-5}$ & - & 1.82 $\cdot 10^{-3}$ & - \\ 
32 $\times$ 32 & 8.35 $\cdot 10^{-6}$ & 3.36 & 4.00 $\cdot 10^{-4}$ & 2.19 \\ 
64 $\times$ 64 & 8.09 $\cdot 10^{-7}$ & 3.37 & 1.16 $\cdot 10^{-5}$ & 5.11 \\ 
128 $\times$ 128 & 1.23 $\cdot 10^{-7}$ & 2.71 & 1.61 $\cdot 10^{-6}$ & 2.85 \\ 
256 $\times$ 256 & 5.32 $\cdot 10^{-8}$ & 1.22 & 5.62 $\cdot 10^{-7}$ & 1.52 \\ 
\hline \hline 
No.~of points & $L^1$ error of $|\nabla c|$ & order & $L^\infty$ error of $|\nabla c|$ & order \\ 
\hline 
16 $\times$ 16 & 7.53 $\cdot 10^{-4}$ & - & 3.51 $\cdot 10^{-3}$ & - \\ 
32 $\times$ 32 & 1.81 $\cdot 10^{-4}$ & 2.05 & 1.48 $\cdot 10^{-3}$ & 1.24 \\ 
64 $\times$ 64 & 4.22 $\cdot 10^{-5}$ & 2.11 & 9.41 $\cdot 10^{-5}$ & 3.98 \\ 
128 $\times$ 128 & 1.06 $\cdot 10^{-5}$ & 1.99 & 2.44 $\cdot 10^{-5}$ & 1.94 \\ 
256 $\times$ 256 & 2.68 $\cdot 10^{-6}$ & 1.98 & 6.84 $\cdot 10^{-6}$ & 1.84 \\ 
\hline \hline 
\end{tabular} 
    \label{table:test2cN}  
 \end{table}

\subsection{Multigrid efficiency}\label{ssec:numtest:multigrid}
In this section we solve Eq.~\eqref{maineq} with data \eqref{dataeq} and the following initial and Dirichlet boundary conditions:
\begin{equation}\label{BCnumtests}
\quad s_0(\vec{x}) = 0, \quad c_0(\vec{x}) = 10, \quad s_B = 1.
\end{equation}
We perform two tests: {\sc Test 3} and {\sc Test 4}, with the domain represented by the level-sets \eqref{test13}  and  \eqref{test24}, respectively  (see Fig.~\ref{fig:domains}).

\begin{figure}[!hbt]
   	\centering
   	\captionsetup{width=0.80\textwidth}
		\includegraphics[width=0.99\textwidth]{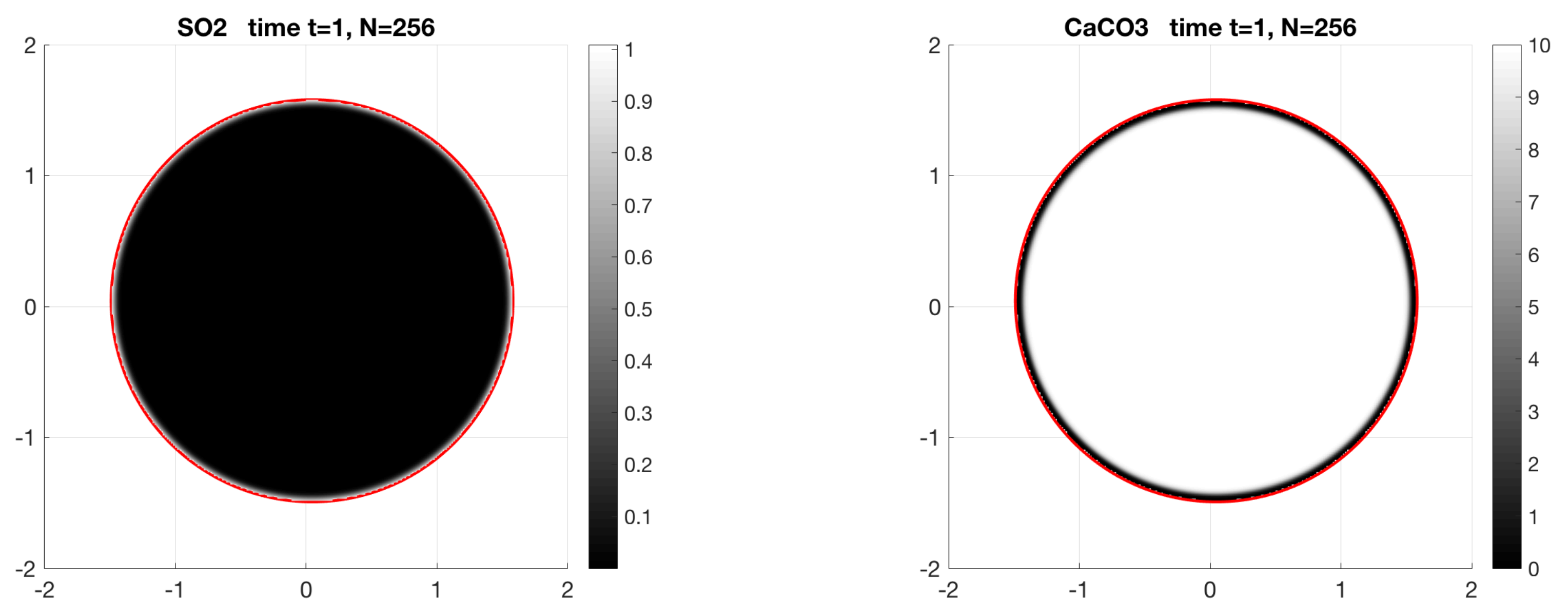}
		\caption{ \footnotesize{Solutions for {\sc Test 3} at time $t=1$ with $N=256$.}}
	\label{fig:solutionTest3}
\end{figure}

\begin{figure}[!hbt]
   	\centering
   	\captionsetup{width=0.80\textwidth}
		\includegraphics[width=0.99\textwidth]{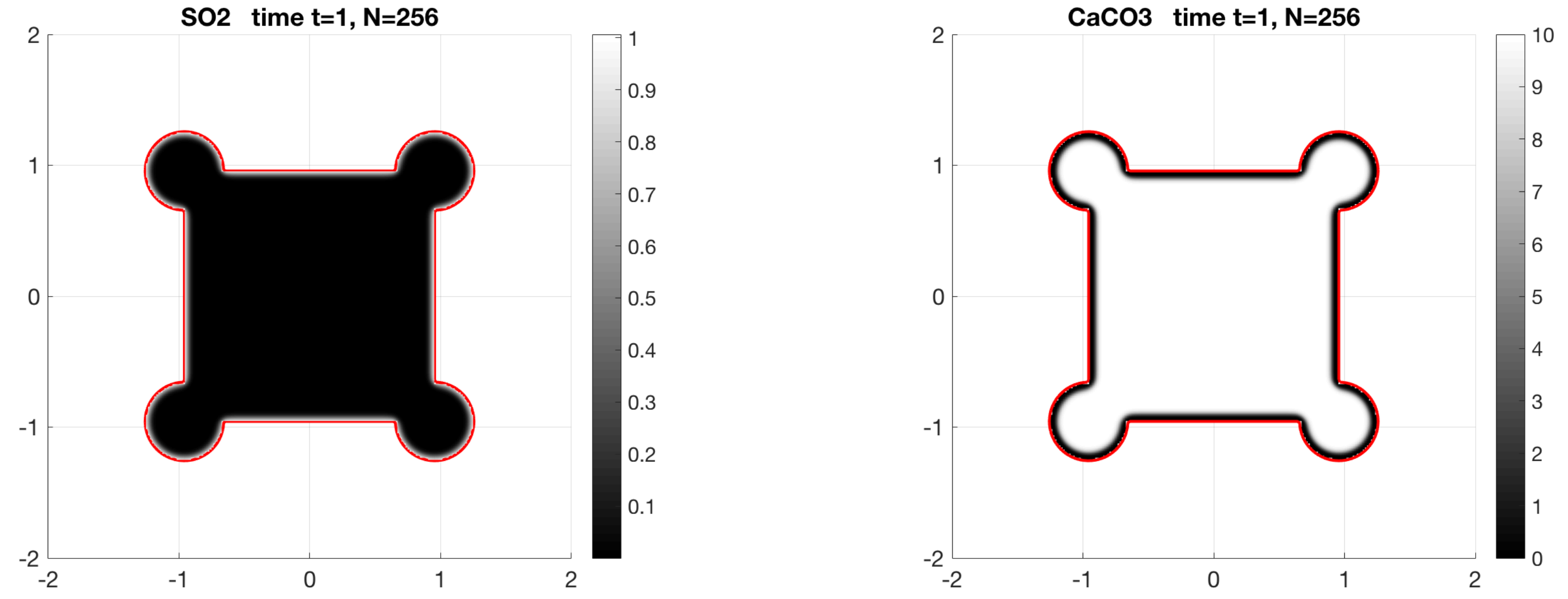}
		\caption{ \footnotesize{Solutions for {\sc Test 4} at time $t=1$ with $N=256$.}}
	\label{fig:solutionTest4}
\end{figure}

\begin{figure}[!hbt]
 \begin{minipage}[c]{0.99\textwidth}
   	\centering
   	\captionsetup{width=0.80\textwidth}
		\includegraphics[width=0.69\textwidth]{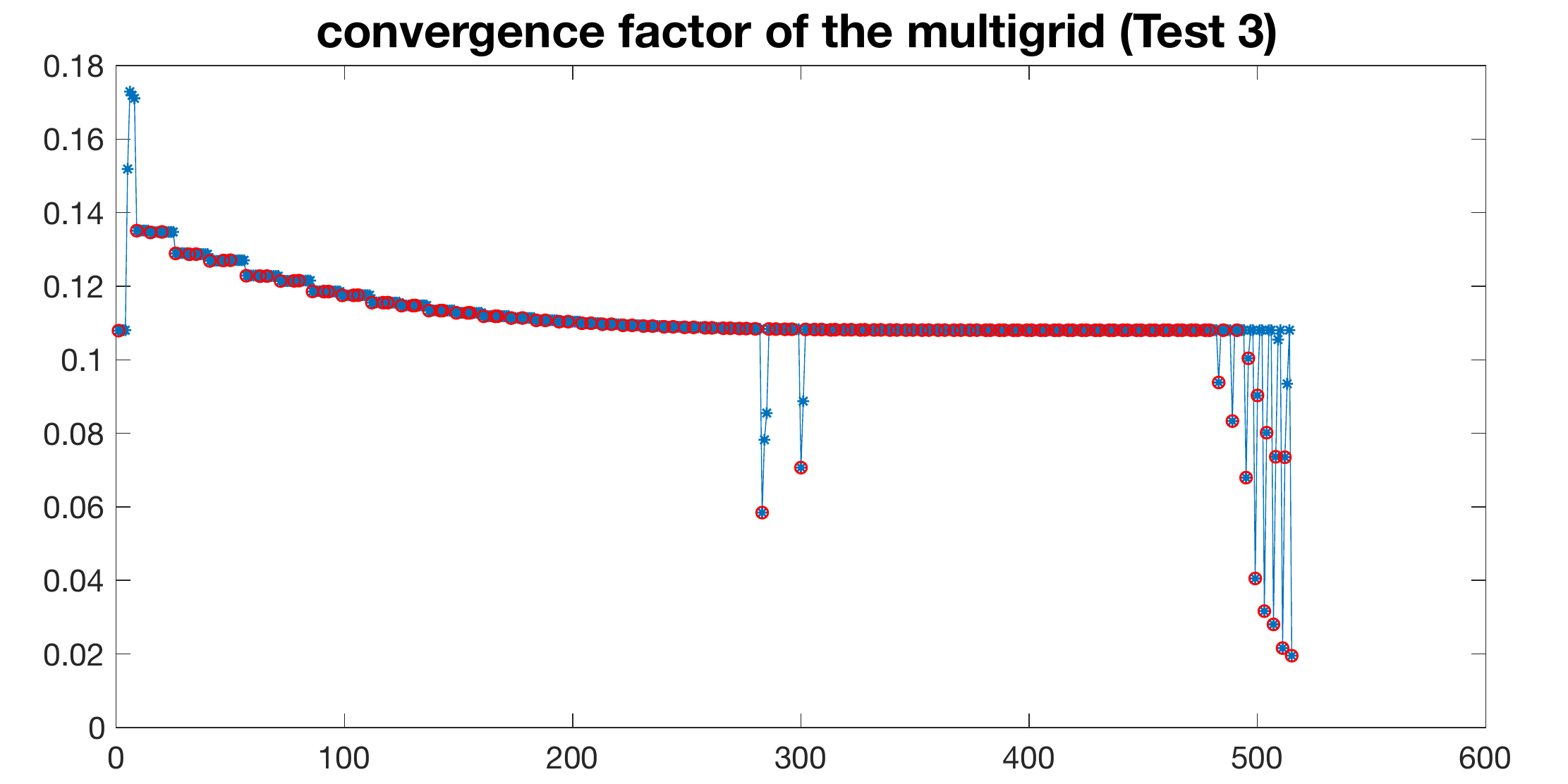} 
 \end{minipage}
 \ \hspace{2mm} \hspace{3mm} \
 \begin{minipage}[c]{0.99\textwidth}
   	\centering
   	\captionsetup{width=0.80\textwidth}
		\includegraphics[width=0.69\textwidth]{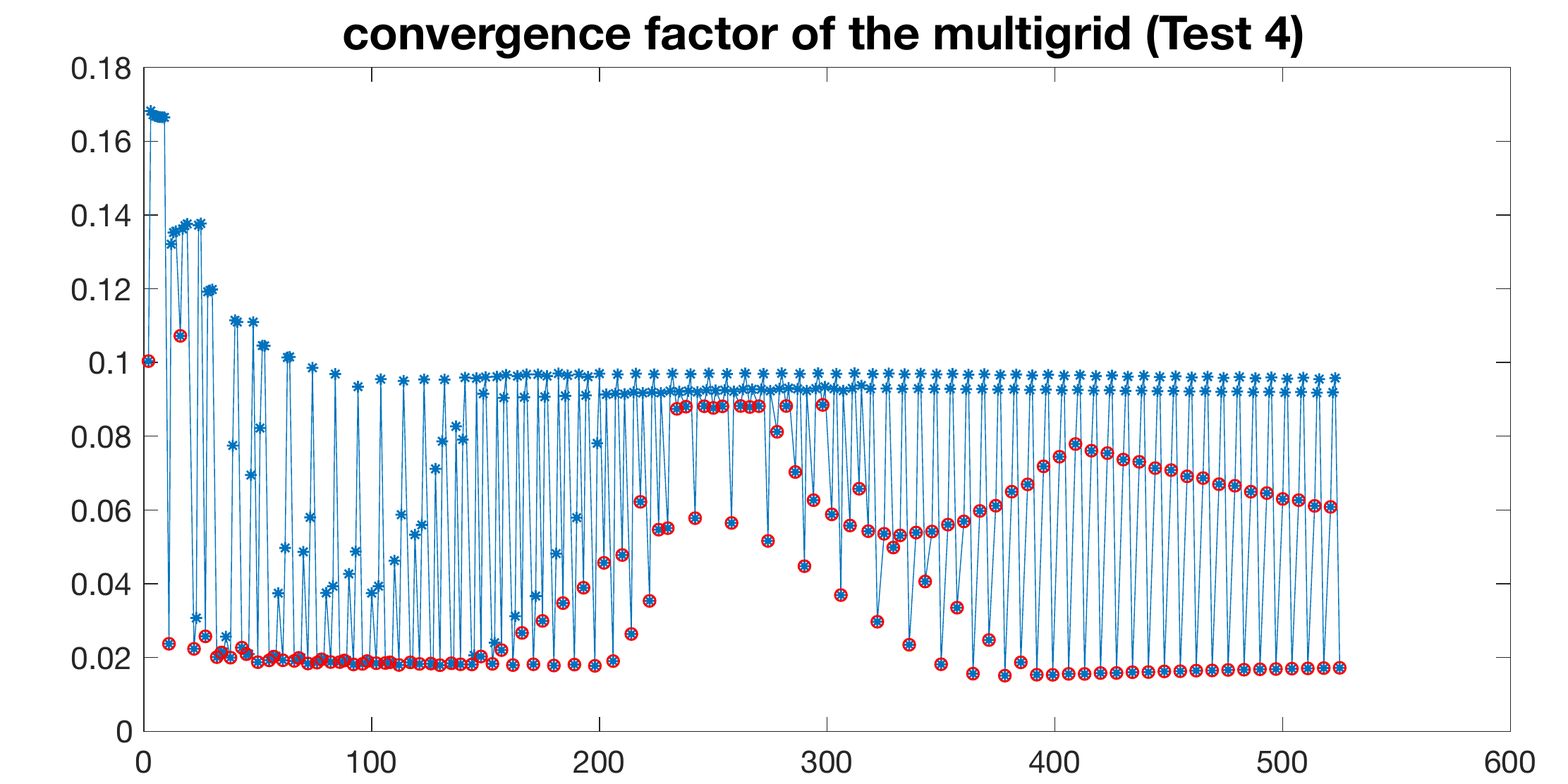} 
 \end{minipage}
 		\caption{ \footnotesize{Convergence factors versus the W-cycle iterations for {\sc Test 3} (top) and {\sc Test 4} (bottom). The first convergence factor obtained in each linear system \eqref{NewtonLS} (either of the same Newton-Raphson step or a new time step) is circled (in red). The convergence factors of the first few linear systems are slightly higher due to the inconsistency of the initial and boundary conditions for $s$. After a few linear systems, the convergence factors are mainly distributed around $\rho=0.119$, which is the predicted value by the Local Fourier Analysis for scalar multigrid in rectangular domains, showing that the multigrid efficiency has not been degraded by the non-rectangular domain and the ghost-point approach.} }
 	\label{fig:rhohist}
\end{figure}

Solutions at time $t=1$ are plotted in Figs.~\ref{fig:solutionTest3} and \ref{fig:solutionTest4}. For each W-cycle of the multigrid method, we compute the convergence factor as:
\[
\rho_\text{MG}^{(q)} = \frac{\left\| \vec{r}_h^{(q)}    \right\|_\infty}{\left\| \vec{r}_h^{(q-1)}   \right\|_\infty},
\]
where 
\begin{equation*}
\vec{r}_h^{(q)} =
\begin{bmatrix}
\vec{r}^s_h \\
\vec{r}^c_h 
\end{bmatrix}
=
\begin{bmatrix}
F^s (W^{(n+1,k)})- \left( J^{ss} (W^{(n+1,k)}) \cdot \Delta \vec{s} + J^{sc}(W^{(n+1,k)}) \cdot \Delta \vec{c} \right) \\
F^c (W^{(n+1,k)})- \left( J^{cs} (W^{(n+1,k)}) \cdot \Delta \vec{s} + J^{cc}(W^{(n+1,k)}) \cdot \Delta \vec{c} \right)
\end{bmatrix}
\end{equation*}
is the defect after $q$ W-cycles.
Convergence factors are plotted in Fig.~\ref{fig:rhohist} versus the W-cycle iterations. The first convergence factor obtained in each linear system \eqref{NewtonLS} (either of the same Newton-Raphson step or a new time step) is circled (in red). The convergence factors of the first few linear systems are slightly higher due to the inconsistency of the initial and boundary conditions for $s$. After a few linear systems, the convergence factors are mainly distributed around $\rho=0.119$, which is the predicted value by the Local Fourier Analysis for scalar multigrid in rectangular domains \cite[Table 4.1, page 117]{Trottemberg:MG}, showing that the multigrid efficiency has not been degraded by the non-rectangular domain and the ghost-point approach.


\subsection{Complex geometries}\label{ssec:numtest:realgeometry}

\begin{figure}[!hbt]
   	\centering
   	\captionsetup{width=0.80\textwidth}
		\includegraphics[width=0.99\textwidth]{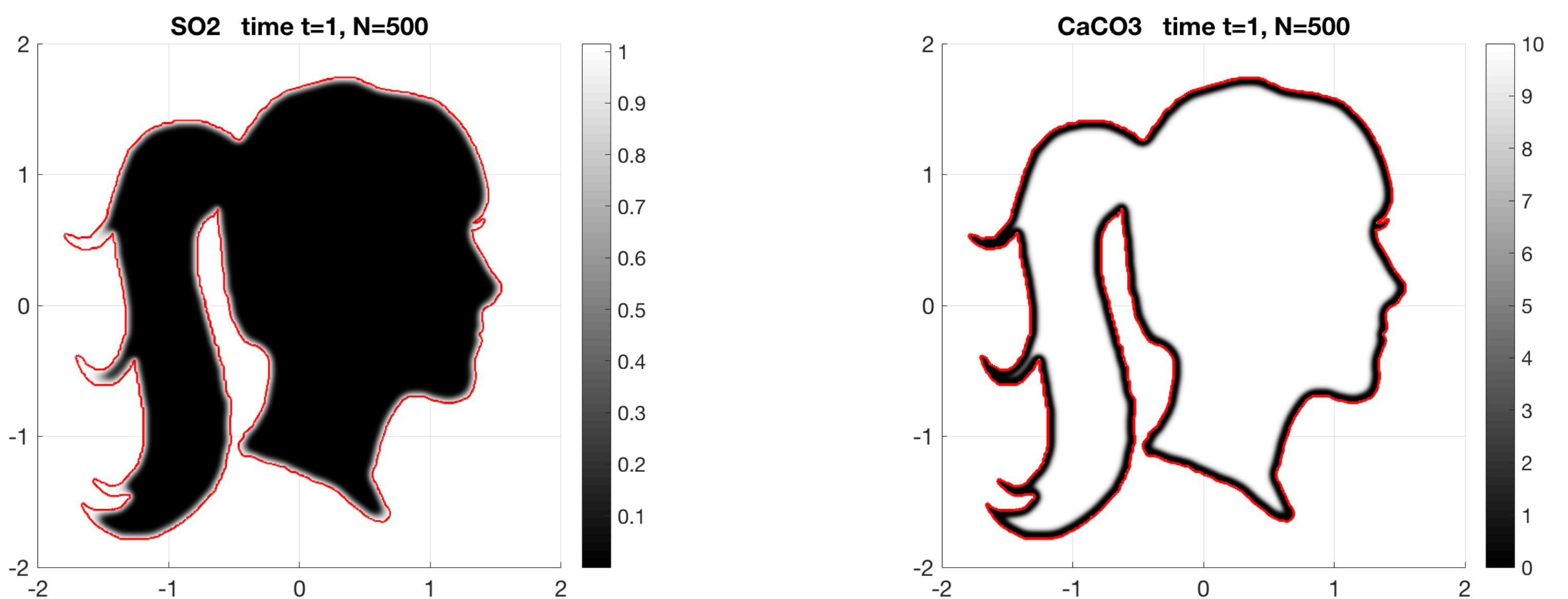}
		\caption{ \footnotesize{Solutions for {\sc Test 5} at time $t=1$ with $N=512$.} }
	\label{fig:solutionTest5}
\end{figure}

\begin{figure}[!hbt]
   	\centering
   	\captionsetup{width=0.80\textwidth}
		\includegraphics[width=0.99\textwidth]{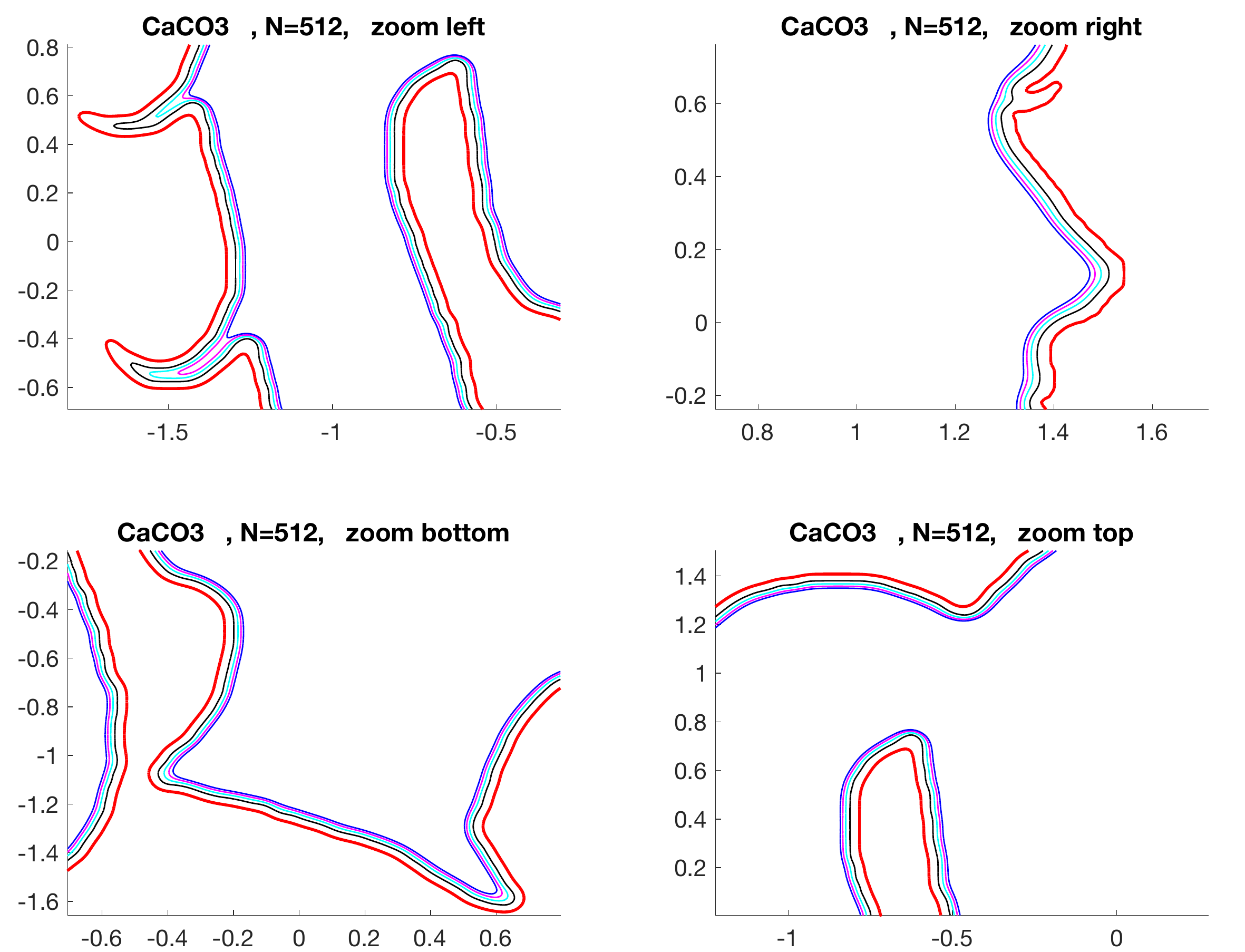}
		\caption{ \footnotesize{Contour plot of $c=5$ ($\text{CaCO}_3$) for {\sc Test 5} at times $t=0.25$, $t=0.5$, $t=0.75$ and $t=1$. The boundary of the domain is the most external line, while the four contour lines can be easily recognised as they move away from the boundary as $t$ increases. In the electronic version of the paper, they are colored black ($t=0.25$), cyan ($t=0.5$), magenta ($t=0.75$) and blue ($t=1$). The boundary of the domain is red. }}
	\label{fig:zoomTest5}
\end{figure}

\begin{figure}[!hbt]
   	\centering
   	\captionsetup{width=0.80\textwidth}
		\includegraphics[width=0.99\textwidth]{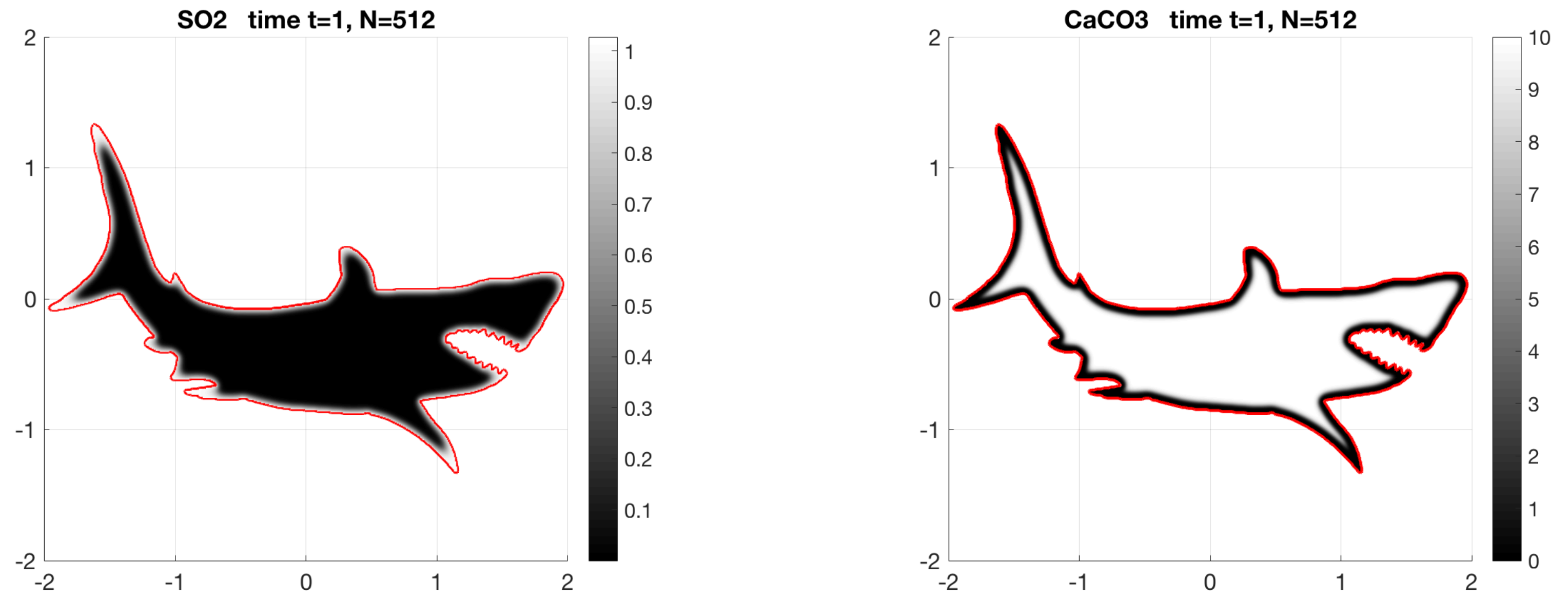}
		\caption{ \footnotesize{Solutions for {\sc Test 6} at time $t=1$ with $N=512$.} }
	\label{fig:solutionTest6}
\end{figure}

\begin{figure}[!hbt]
   	\centering
   	\captionsetup{width=0.80\textwidth}
		\includegraphics[width=0.99\textwidth]{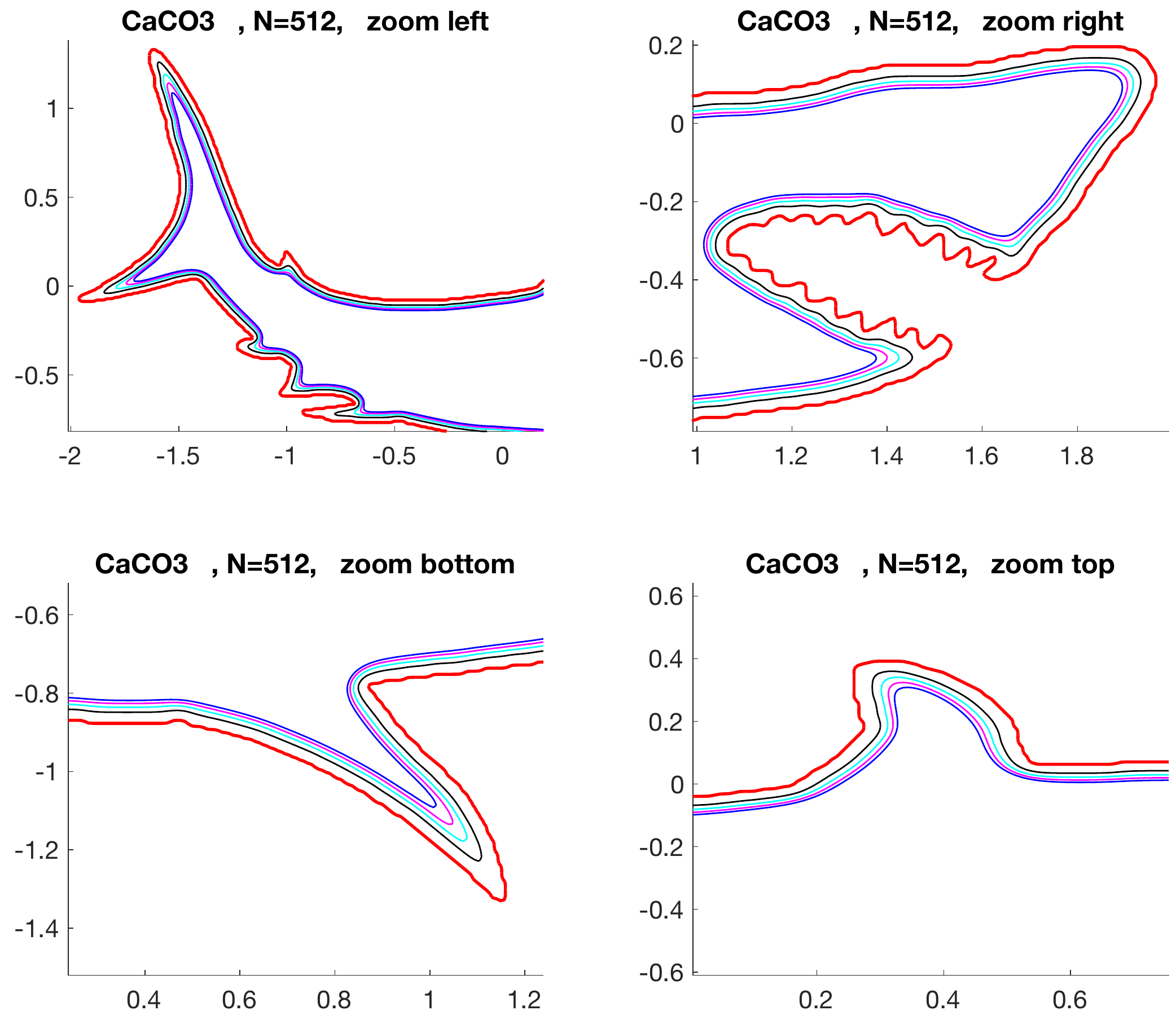}
		\caption{ \footnotesize{Contour plot of $c=5$ ($\text{CaCO}_3$) for {\sc Test 6} at times $t=0.25$, $t=0.5$, $t=0.75$ and $t=1$. The boundary of the domain is the most external line, while the four contour lines can be easily recognised as they move away from the boundary as $t$ increases. In the electronic version of the paper, they are colored black ($t=0.25$), cyan ($t=0.5$), magenta ($t=0.75$) and blue ($t=1$). The boundary of the domain is red. }}
	\label{fig:zoomTest6}
\end{figure}

\begin{figure}[!hbt]
   	\centering
   	\captionsetup{width=0.80\textwidth}
			\includegraphics[width=0.99\textwidth]{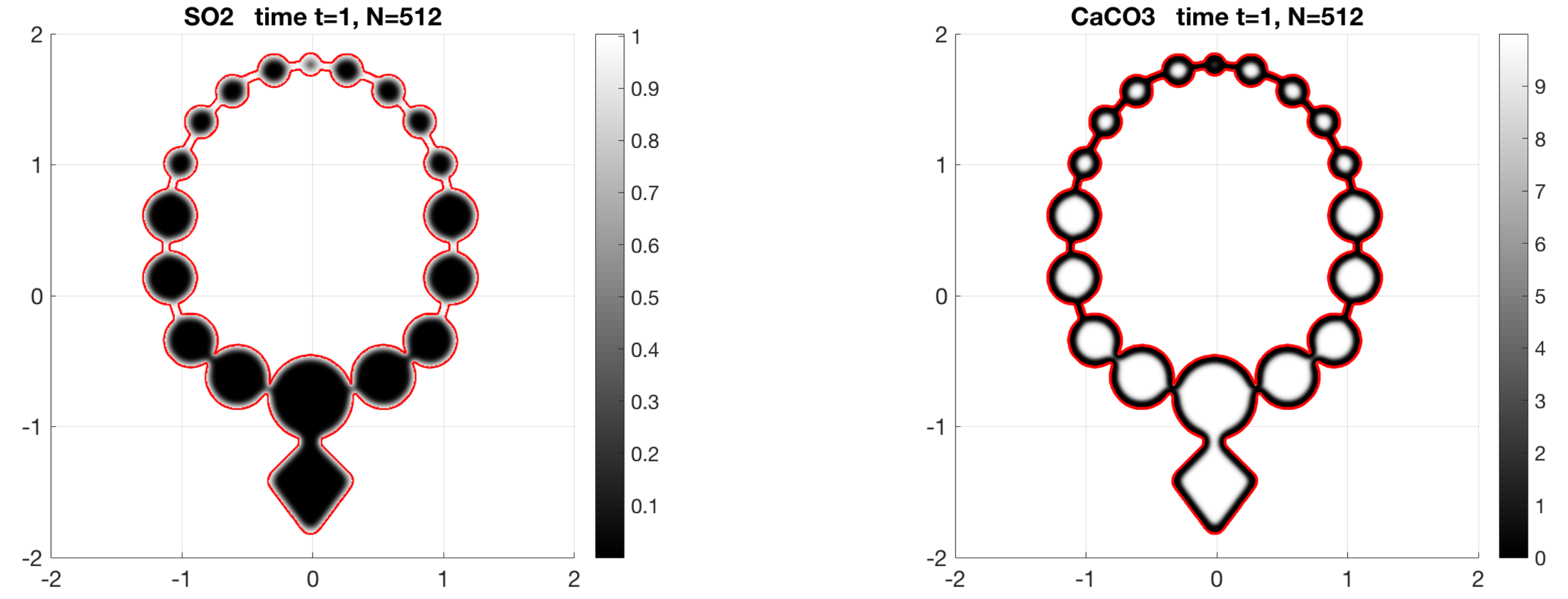}
		\caption{ \footnotesize{Solutions for {\sc Test 7} at time $t=1$ with $N=512$.} }
	\label{fig:solutionTest7}
\end{figure}

\begin{figure}[!hbt]
   	\centering
   	\captionsetup{width=0.80\textwidth}
		\includegraphics[width=0.99\textwidth]{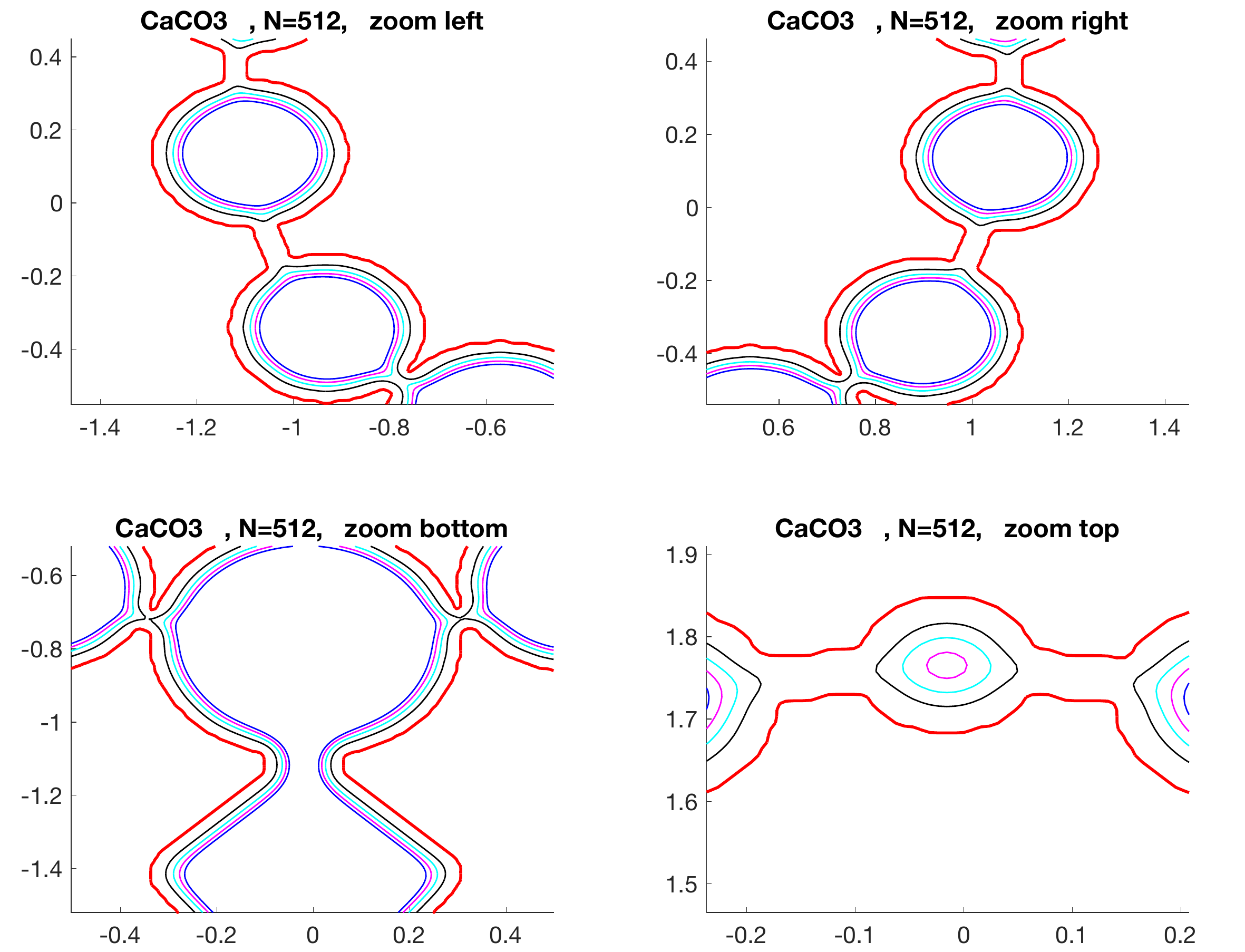}
		\caption{ \footnotesize{Contour plot of $c=5$ ($\text{CaCO}_3$) for {\sc Test 7} at times $t=0.25$, $t=0.5$, $t=0.75$ and $t=1$. The boundary of the domain is the most external line, while the four contour lines can be easily recognised as they move away from the boundary as $t$ increases. In the electronic version of the paper, they are colored black ($t=0.25$), cyan ($t=0.5$), magenta ($t=0.75$) and blue ($t=1$). The boundary of the domain is red. }}
	\label{fig:zoomTest7}
\end{figure}

In this section we show how the method performs on more complex geometries.
We use data \eqref{dataeq} and \eqref{BCnumtests} and we test three domains: a woman head profile ({\sc Test 5}), a shark ({\sc Test 6}) and a necklace ({\sc Test 7}). The solutions at time $t=1$ are represented in Figs.~\ref{fig:solutionTest5}, \ref{fig:solutionTest6} and \ref{fig:solutionTest7}, respectively, while zooms on some relevant regions, with the contour plot of $\text{CaCO}_3$ at times $t=0.25$, $t=0.5$, $t=0.75$ and $t=1$, are represented in Figs.~\ref{fig:zoomTest5}, \ref{fig:zoomTest6} and \ref{fig:zoomTest7}, respectively. In general, we can see that the reaction is quicker around {\it corner points}, i.e.~regions where the boundary has a higher curvature. For example, in {\sc Test 5} the eyelash of the woman is entirely gypsum already at time $t=0.25$ (top-right plot of Fig.~\ref{fig:zoomTest5}), while the hair strands progressively draw back at times $t=0.25$, $t=0.5$ and $t=0.75$, and almost disappear at time $t=1$ (top-left plot of Fig.~\ref{fig:zoomTest5}). The gypsum formation is also quicker around reentrant corners, as we can see in the bottom-left and bottom-right plots of Fig.~\ref{fig:zoomTest5}, for example.
Similar conclusions may be drawn for the shark in {\sc Test 6}, where the teeth (top-right plot of Fig.~\ref{fig:zoomTest6}) and the secondary dorsal fin (top-left plot of Fig.~\ref{fig:zoomTest6}) are transformed quickly. The damage on the caudal fin is quicker around the tips (top-left plot of Fig.~\ref{fig:zoomTest6}) and around the reentrant corner of the pectoral fin (bottom-left plot of Fig.~\ref{fig:zoomTest6}) and of the primary dorsal fin (bottom-right plot of Fig.~\ref{fig:zoomTest6}).
{\sc Tests 5} and {\sc 6} confirm the qualitative observation made by practitioners in the field of conservation of cultural heritage that the parts of a manufact that are most quickly affected by the sulfation of marble are the higher details of the decoration and the sharp edges. It is important to take into account that gypsum is soluble in water and very prone to breaking due to thermal shocks and thus that the areas with high gypsum content, in a real case, would be quickly lost by dissolution into rainwater or by exfoliation by the dilatation due to the cyclic seasonal temperature variations.
{\sc Test 7} shows qualitatively analogous results. In fact, the lace connecting the beads is the most quickly damaged part, due to its reduced thickness, followed by the beads in order of size. The larger beads and the diamond shaped pendant are less sulfated and suffer damages almost only close to the reentrant corners. 
Fig.~\ref{fig:solutionTest7} shows that in this case gypsum dissolution in water or its exfoliation would cause a dramatic topological change, disconnecting the necklace into many separated pieces.



\section{Conclusion}\label{sec:conclusion}
 Having in mind the modeling of marble degradation under chemical pollutants e.g.~the sulfation process, we considered the governing nonlinear equations and their numerical approximation. The space domain is implicitly defined using a level-set approach. We employed a Crank-Nicolson in time, while for the space variables the discretization is performed by a standard Finite-Difference scheme for grid points inside the domain  and by a ghost-cell technique on the ghost points (by using boundary conditions).

The solution of the large nonlinear system has been obtained by a Newton-Raphson procedure and by a tailored multigrid technique. All the numerical experiments have given very satisfactory results both from the viewpoint of the reconstruction quality and of the computational efficiency.

As future steps we can include, from the numerical analysis point of view, the spectral analysis of the resulting matrices from a GLT viewpoint \cite{gltbook1,gltbook2} having in mind a rigorous convergence analysis of the considered multigrid techniques.

From a modelling point of view, it would be interesting to extend the computational techniques introduced in this paper to the models of degradation processes that employ an evolving domain (e.g. \cite{CFN:08:swelling}) or that include internal moving interfaces among layers of materials with different properties (e.g. \cite{CdFN:14:copper,Nikolopoulos:14:mushy}), to take into account some effects like swelling and corrosion.
In this respect, we point out that the level-set technique introduced in this paper would be able to track correctly the pristine marble domain even if it was disconnected during the time evolution as in the example of Fig.~\ref{fig:solutionTest7}.

From a computational point of view, the method will be extended to three dimensional problems, where a realistic piece of work will be modelled from a laser scanner 3D reconstruction. To reduce the computational cost, the computational strategy will include an Adaptive Meshing Refinement (AMR) approach, where the Cartesian mesh will be refined around critical regions such as the boundary and/or the internal moving interface, by extending AMR strategies already developed for simpler problems~\cite{mirzadeh2011second, theillard2013multigrid}.

\section*{Acknowledgments} The work has been partially supported by the London Mathematical Society Computer Science Small Grants -- Scheme 7 (Ref.~SC7-1617-02) and the Research in Pairs -- Scheme 4 (Ref.~41739).

\addcontentsline{toc}{chapter}{References}
\bibliographystyle{jabbrv_abbrv}
\bibliography{bibliography}
\end{document}